\documentclass[11pt]{article}
\usepackage{enumerate}
\usepackage{amssymb,a4wide,latexsym,makeidx,epsfig,fleqn}
\usepackage{amsthm}
\usepackage{amsmath}
\usepackage{enumerate}

\newtheorem{Exa}{Example}[section]

\begin{document}
\textwidth 150mm \textheight 225mm
\title{A new \emph{S}-type eigenvalue localization set for tensors and its applications
\thanks{ Supported by the National Natural Science Foundation of China~(No.11171273).}}
\author{{Zhengge Huang, Ligong Wang\footnote{Corresponding author.} , Zhong Xu and Jingjing Cui}\\
{\small Department of Applied Mathematics, School of Science, Northwestern
Polytechnical University,}\\ {\small  Xi'an, Shaanxi 710072,
People's Republic
of China.}\\
{\small E-mails: ZhenggeHuang@mail.nwpu.edu.cn; lgwang@nwpu.edu.cn(or lgwangmath@163.com);}\\
{\small zhongxu@nwpu.edu.cn; JingjingCui@mail.nwpu.edu.cn}\\
 }
\date{}
\maketitle
\begin{center}
\begin{minipage}{120mm}
\vskip 0.3cm
\begin{center}
{\small {\bf Abstract}}
\end{center}
{\small A new \emph{S}-type eigenvalue localization set for tensors is derived by breaking $N=\{1,2,\cdots,n\}$ into disjoint subsets $S$ and its complement. It is proved that this new set is tighter than those presented by Qi (Journal of Symbolic Computation
40 (2005) 1302-1324), Li et al. (Numer. Linear Algebra Appl. 21 (2014) 39-50) and Li et al. (Linear Algebra Appl. 493 (2016) 469-483).  As applications, checkable sufficient conditions for the positive definiteness and the positive semi-definiteness of tensors are proposed. Moreover, based on this new set, we establish a new upper bound for the spectral radius of nonnegative tensors and a lower bound for the minimum \emph{H}-eigenvalue of weakly irreducible strong \emph{M}-tensors in this paper. We demonstrate that these bounds are sharper than those obtained by Li et al. (Numer.
Linear Algebra Appl. 21 (2014) 39-50) and He and Huang (J. Inequal. Appl.
114 (2014) 2014). Numerical examples are also given to illustrate this fact.

 }
\end{minipage}
\end{center}

\section{Introduction }
\label{sec:ch6-introduction}

Eigenvalue problems of higher order tensors have become an
important topic in applied mathematics branch, numerical multilinear algebra, and it
has a wide range of practical applications, such as best-rank one approximation in data analysis \cite{23}, higher order Markov chains \cite{24}, molecular conformation \cite{25} and so forth. Recently, tensor eigenvalues have received much attention in the literatures \cite{1,2,14,3,18,17,4,5,6,10,7,8,9,11,12,13}.

One of many practical applications of eigenvalues of tensors is that one can use the smallest \emph{H}-eigenvalue of an even-order real symmetric tensor to identify its positive
(semi-)definiteness, consequently, can identify the positive (semi-)definiteness of the multivariate homogeneous polynomial determined by this tensor, for details, see \cite{1,21,32}.

However, as mentioned in \cite{16,21,33}, it is not easy to compute the smallest \emph{H}-eigenvalue of tensors when the order and dimension are very large, we always try to give a set including all eigenvalues in
the complex. Some sets including all eigenvalues of tensors have been presented by some researchers \cite{1,13,21,32,16,33,22,31}. In particular, if one of these sets for an even-order real symmetric tensor is in the
right-half complex plane, then we can conclude that the smallest \emph{H}-eigenvalue is positive,
consequently, the corresponding tensor is positive definite. Therefore, the main aim of this paper is to study the new eigenvalue localization set for tensors called new \emph{S}-type eigenvalue localization set, which is sharper than some existing ones.

For a positive integer $n,N$ denotes the set $N=\{1,2,\cdots,n\}$. The set of all real numbers is denoted by $\mathbb{R}$, and $\mathbb{C}$ denotes the set of all complex numbers. Here, we call $\mathcal{A}=(a_{i_{1}\cdots i_{m}})$ a complex (real) tensor of order $m$ dimension $n$, denoted by $\mathbb{C}^{[m,n]}$$(\mathbb{R}^{[m,n]})$, if $a_{i_{1}\cdots i_{m}}\in{\mathbb{C}}(\mathbb{R})$, where $i_{j}\in{N}$ for $j=1,2,\cdots,m$ \cite{16}.

Let $\mathcal{A}\in\mathbb{R}^{[m,n]}$, and $x\in{\mathbb{C}}^{n}$. Then
\begin{eqnarray}
\mathcal{A}x^{m-1}:=\bigg(\sum\limits_{i_{2},\cdots, i_{m}=1}^{n}a_{ii_{2}\cdots i_{m}}x_{i_{2}}\cdots x_{i_{m}}\bigg)_{1\leq{i}\leq{n}},\nonumber
\end{eqnarray}
a pair $(\lambda,x)\in{\mathbb{C}}\times(\mathbb{C}^{n}/\{0\})$ is called an eigenpair of $\mathcal{A}$ \cite{9} if
\begin{eqnarray}
\mathcal{A}x^{m-1}=\lambda x^{[m-1]},\nonumber
\end{eqnarray}
where $x^{[m-1]}=(x_{1}^{m-1},x_{2}^{m-1},\cdots,x_{n}^{m-1})^{T}$ \cite{22}. Furthermore, we call $(\lambda,x)$ an $H$-eigenpair, if both $\lambda$ and $x$ are real \cite{1}.

A real tensor of order $m$ dimension $n$ is called the unit tensor \cite{21}, denoted by $\mathcal{I}$, if its entries are $\delta_{i_{1}\cdots i_{m}}$ for $i_{1},\cdots, i_{m}\in{N}$, where
\begin{equation}
  \delta_{i_{1}\cdots i_{m}}=\left\{
   \begin{aligned}
   &1,\quad \mathrm{if}\ i_{1}=\cdots=i_{m},\\
   &0,\quad \mathrm{otherwise}. \\
   \end{aligned}
   \right.\nonumber
\end{equation}

An $m$-order $n$-dimensional tensor $\mathcal{A}$ is called nonnegative \cite{3,4,5,14,20}, if each entry is nonnegative. We call a tensor $\mathcal{A}$ as a $\mathrm{Z}$-tensor, if all of its off-diagonal entries are non-positive, which is equivalent to write $\mathcal{A}=s\mathcal{I}-\mathcal{B}$, where $s>0$ and $\mathcal{B}$ is a nonnegative tensor ($\mathcal{B}\geq{0}$), denote by $\mathbb{Z}$ the set of $m$-order and $n$-dimensional $\mathrm{Z}$-tensors. A $\mathrm{Z}$-tensor $\mathcal{A}=s\mathcal{I}-\mathcal{B}$ is an \emph{M}-tensor if $s\geq\rho(\mathcal{B})$, and it is a nonsingular (strong) \emph{M}-tensor if $s>\rho(\mathcal{B})$ \cite{12,19,34}.

A tensor $\mathcal{A}=(a_{i_{1}\cdots i_{m}})\in{\mathbb{C}}^{[m,n]}$ is called weakly reducible, if there exists a nonempty proper index subset $I\subset N$ such that
$a_{i_{1}i_{2}\cdots i_{m}} = 0$, $\forall i_{1}\in I$, $\exists i_{j}\notin{I}$, $j=2,\cdots,n$.
If $\mathcal{A}$ is not weakly reducible, then we call $\mathcal{A}$ weakly irreducible \cite{5,13}. The tensor $\mathcal{A}$ is called reducible if there exists a nonempty proper index subset $\mathbb{J}\subset N$ such that $a_{i_{1}i_{2}\cdots i_{m}}=0$, $\forall i_{1}\in\mathbb{J}$, $\forall i_{2},\cdots,i_{m}\notin\mathbb{J}$. If $\mathcal{A}$ is not reducible, then we call $\mathcal{A}$ is irreducible \cite{11}.
The spectral radius $\rho(\mathcal{A})$ \cite{5} of the tensor $\mathcal{A}$ is defined as
\begin{eqnarray}
\rho(\mathcal{A})=\max\{|\lambda|:\lambda\ \mathrm{is}\ \mathrm{an}\ \mathrm{eigenvalue}\ \mathrm{of}\ \mathcal{A}\}.\nonumber
\end{eqnarray}
Denoted by $\tau(\mathcal{A})$ the minimum value of the real part of all eigenvalues of the strong \emph{M}-tensor $\mathcal{A}$ \cite{6}. A real tensor $\mathcal{A}=(a_{i_{1}\cdots i_{m}})$ is called symmetric \cite{1,4,13,16,31,32} if
\begin{eqnarray}
a_{i_{1}\cdots i_{m}}=a_{\pi(i_{1}\cdots i_{m})},\ \forall \pi\in\Pi_{m},\nonumber
\end{eqnarray}
where $\Pi_{m}$ is the permutation group of $m$ indices.

Let $\mathcal{A}=(a_{i_{1}\cdots i_{m}})\in{\mathbb{R}}^{[m,n]}$. For $i,j\in{N}$, $j\neq{i}$, denote
\begin{eqnarray}
&&R_{i}(\mathcal{A})=\sum\limits_{i_{2},\cdots, i_{m}=1}^{n}a_{ii_{2}\cdots i_{m}},\ R_{\max}(\mathcal{A})=\max\limits_{i\in{N}}R_{i}(\mathcal{A}),\ R_{\min}(\mathcal{A})=\min\limits_{i\in{N}}R_{i}(\mathcal{A}),\nonumber\\
&&r_{i}(\mathcal{A})=\sum\limits_{\delta_{ii_{2}\cdots i_{m}}=0}|a_{ii_{2}\cdots i_{m}}|,\
r_{i}^{j}(\mathcal{A})=\sum\limits_{\substack{\delta_{ii_{2}\cdots i_{m}}=0,\\ \delta_{ji_{2}\cdots i_{m}}=0}}|a_{ii_{2}\cdots i_{m}}|=r_{i}(\mathcal{A})-|a_{ij\cdots{j}}|.\nonumber
\end{eqnarray}

Recently, many literatures have been focused on the bounds of the spectral radius
of nonnegative tensor in \cite{5,10,7,8,9,11,28}. Also, in \cite{6}, He and Huang obtained the upper
and lower bounds for the minimum \emph{H}-eigenvalue of irreducible strong \emph{M}-tensors. Wang and Wei \cite{7} presented some new bounds for the minimum \emph{H}-eigenvalue of weakly irreducible strong \emph{M}-tensors, and showed those are better than the ones in \cite{6} in some cases.
Based on the new set established in this paper, the other main results of this paper is to provide
sharper bounds for the spectral radius of nonnegative tensors and the minimum
\emph{H}-eigenvalue of weakly irreducible nonsingular \emph{M}-tensors, which improve some existing ones.

Before presenting our results, we review the existing results related to the eigenvalue localization sets for tensors. In 2005, Qi \cite{6} generalized $\mathrm{Ger\check{s}gorin}$ eigenvalue localization theorem from matrices to real supersymmetric tensors, which can be easily extended to general tensors \cite{4,31}.
\newtheorem{lem1}{Lemma}[section]
\begin{lem1}\emph{\cite{1}}\label{lem1}
Let $\mathcal{A}=(a_{i_{1}\cdots i_{m}})\in{\mathbb{C}}^{[m,n]}$, $n\geq{2}$. Then
\begin{eqnarray}
\sigma(\mathcal{A})\subseteq\Gamma(\mathcal{A})=\bigcup\limits_{i\in{N}}\Gamma_{i}(\mathcal{A}),\nonumber
\end{eqnarray}
where $\sigma(\mathcal{A})$ is the set of all the eigenvalues of $\mathcal{A}$ and
\begin{eqnarray}
\Gamma_{i}(\mathcal{A})=\{z\in{\mathbb{C}}:|z-a_{i\cdots i}|\leq r_{i}(\mathcal{A})\}.\nonumber
\end{eqnarray}
\end{lem1}
To get sharper eigenvalue localization sets than $\Gamma(\mathcal{A})$, Li et al. \cite{31} extended the Brauer's eigenvalue localization set of matrices \cite{26} and proposed the following Brauer-type eigenvalue localization set for tensors.
\newtheorem{lem2}[lem1]{Lemma}
\begin{lem2}\emph{\cite{31}}\label{lem1}
Let $\mathcal{A}=(a_{i_{1}\cdots i_{m}})\in{\mathbb{C}}^{[m,n]}$, $n\geq{2}$. Then
\begin{eqnarray}
\sigma(\mathcal{A})\subseteq\mathcal{K}(\mathcal{A})=\bigcup\limits_{i,j\in{N},j\neq{i}}\mathcal{K}_{i,j}(\mathcal{A}),\nonumber
\end{eqnarray}
where
\begin{eqnarray}
\mathcal{K}_{i,j}(\mathcal{A})=\{z\in{\mathbb{C}}:(|z-a_{i\cdots i}|-r_{i}^{j}(\mathcal{A}))|z-a_{j\cdots j}|\leq |a_{ij\cdots j}|r_{j}(\mathcal{A})\}.\nonumber
\end{eqnarray}
\end{lem2}
In addition, in order to reduce computations of determining the sets $\sigma(\mathcal{A})$, Li et al. \cite{31} also presented the following $S$-type eigenvalue localization set by breaking $N$ into disjoint subsets $S$ and $\bar{S}$, where $\bar{S}$ is the complement of $S$ in $N$.
\newtheorem{lem3}[lem1]{Lemma}
\begin{lem3}\emph{\cite{31}}\label{lem1}
Let $\mathcal{A}=(a_{i_{1}\cdots i_{m}})\in{\mathbb{C}}^{[m,n]}$, $n\geq{2}$, and $S$ be a nonempty proper subset of $N$. Then
\begin{eqnarray}
\sigma(\mathcal{A})\subseteq\mathcal{K}^{S}(\mathcal{A})=\left(\bigcup\limits_{i\in{S},j\in{\bar{S}}}\mathcal{K}_{i,j}(\mathcal{A})\right)\bigcup
\left(\bigcup\limits_{i\in{\bar{S}},j\in{{S}}}\mathcal{K}_{i,j}(\mathcal{A})\right),\nonumber
\end{eqnarray}
where $\mathcal{K}_{i,j}(\mathcal{A})$ $(i\in{S},j\in{\bar{S}}\ \mathrm{or}\ i\in{\bar{S}},j\in{{S}})$ is defined as in Lemma 1.2.
\end{lem3}
Very recently, by the technique in \cite{31}, Li et al. \cite{16} gave the new eigenvalue localization set involved with a proper subset $S$ of $N$, and by the following three sets:
\begin{eqnarray}
&&\Delta^{N}=\{(i_{2},i_{3},\cdots,i_{m}):\mathrm{each}\ i_{j}\in{N}\ \mathrm{for}\ j=2,3,\cdots,m\},\nonumber\\
&&\Delta^{S}=\{(i_{2},i_{3},\cdots,i_{m}):\mathrm{each}\ i_{j}\in{S}\ \mathrm{for}\ j=2,3,\cdots,m\},\ \overline{\Delta^{S}}=\Delta^{N}\backslash\Delta^{S}.\nonumber
\end{eqnarray}
\newtheorem{lem10}[lem1]{Lemma}
\begin{lem10}\emph{\cite{16}}\label{lem10}
Let $\mathcal{A}=(a_{i_{1}\cdots i_{m}})\in{\mathbb{C}}^{[m,n]}$, $n\geq{2}$, and $S$ be a nonempty proper subset of $N$. Then
\begin{eqnarray}
\sigma(\mathcal{A})\subseteq\Omega^{S}(\mathcal{A})=\left(\bigcup\limits_{i\in{S},j\in{\bar{S}}}\Omega_{i,j}^{S}(\mathcal{A})\right)\bigcup
\left(\bigcup\limits_{i\in{\bar{S}},j\in{{S}}}\Omega_{i,j}^{\bar{S}}(\mathcal{A})\right),\nonumber
\end{eqnarray}
where
\begin{eqnarray}
&&\Omega_{i,j}^{S}(\mathcal{A})=\{z\in{\mathbb{C}}:(|z-a_{i\cdots i}|)(|z-a_{j\cdots j}|-r_{j}^{\overline{\Delta^{S}}}(\mathcal{A}))\leq r_{i}(\mathcal{A})r_{j}^{\Delta^{S}}(\mathcal{A})\},\nonumber\\
&&\Omega_{i,j}^{\bar{S}}(\mathcal{A})=\{z\in{\mathbb{C}}:(|z-a_{i\cdots i}|)(|z-a_{j\cdots j}|-r_{j}^{\overline{\Delta^{\bar{S}}}}(\mathcal{A}))\leq r_{i}(\mathcal{A})r_{j}^{\Delta^{\bar{S}}}(\mathcal{A})\},\nonumber
\end{eqnarray}
and for $i\in{S}$,
\begin{eqnarray}
r_{i}(\mathcal{A})=r_{i}^{\Delta^{S}}(\mathcal{A})+r_{i}^{\overline{\Delta^{S}}}(\mathcal{A}),\ r_{i}^{j}(\mathcal{A})=r_{i}^{\Delta^{S}}(\mathcal{A})+r_{i}^{\overline{\Delta^{S}}}(\mathcal{A})-|a_{ij\cdots j}|,\nonumber
\end{eqnarray}
with
\begin{eqnarray}
r_{i}^{\Delta^{S}}(\mathcal{A})=\sum\limits_{\substack{(i_{2},\cdots,i_{m})\in{\Delta^{S}},\\ \delta_{ii_{2}\cdots i_{m}=0}}}|a_{ii_{2}\cdots i_{m}}|,\ r_{i}^{\overline{\Delta^{S}}}(\mathcal{A})=\sum\limits_{(i_{2},\cdots,i_{m})\in\overline{\Delta^{S}}}|a_{ii_{2}\cdots i_{m}}|.\nonumber
\end{eqnarray}
\end{lem10}
Theorem 6 in \cite{16} shows that this new set is tighter than the sets $\Gamma(\mathcal{A})$, $\mathcal{K}(\mathcal{A})$ and $\mathcal{K}^{S}(\mathcal{A})$.

In this paper, we focus on investigating the eigenvalue localization sets for tensors, and obtain a new \emph{S}-type eigenvalue localization set for tensors. It is proved to be tighter than the tensor $\mathrm{Ger\check{s}gorin}$ eigenvalue localization set $\Gamma(\mathcal{A})$ in Lemma 1.1, the Brauer's eigenvalue localization set $\mathcal{K}(\mathcal{A})$ in Lemma 1.2, the $S$-type eigenvalue localization set $\mathcal{K}^{S}(\mathcal{A})$ in Lemma 1.3 and another $S$-type eigenvalue localization set $\Omega^{S}(\mathcal{A})$ in Lemma 1.4. As applications, checkable sufficient conditions for the positive definiteness and the positive semi-definiteness of tensors are proposed, and some new bounds for the spectral radius of nonnegative tensors and the minimum \emph{H}-eigenvalue of weakly irreducible strong \emph{M}-tensors are established. The bounds improve some existing ones. Numerical examples are implemented to illustrate this fact.

The outline of this paper is organized as follows. In Section 2, we recollect some useful lemmas
 which are utilized in the next sections. In Section 3, a new \emph{S}-type eigenvalue localization set for tensors is given, and proved to be tighter than the existing ones derived in Lemmas 1.1-1.4. As applications of the results in Section 3, checkable sufficient conditions for the positive definiteness and the positive semi-definiteness of tensors are given in Section 4. Based on the results of Section 3, we propose a new upper bound for the spectral radius of nonnegative tensors in Section 5, comparison results for this new bound and those derived in \cite{31} are also investigated in this section. Section 6 is devoted to exhibit a new lower bound for the minimum \emph{H}-eigenvalue of weakly irreducible strong \emph{M}-tensors, which is proved to be sharper than the ones obtained by He and Huang \cite{6}. Finally, some concluding remarks are given to end this paper in Section 7.

\section{Preliminaries}
\label{sec:Preliminaries}

In this section, we start with some lemmas. They will be useful in the following proofs.
\newtheorem{lem11}[lem1]{Lemma}
\begin{lem11}\emph{\cite{4}}\label{lem11}
If $\mathcal{A}\in{\mathbb{R}}^{[m,n]}$ is nonnegative, then $\rho(\mathcal{A})$ is an
eigenvalue with an entrywise nonnegative eigenvector $x$, i.e., $x\geq0$, $x\neq{0}$, corresponding to it.
\end{lem11}
\newtheorem{lem12}[lem1]{Lemma}
\begin{lem12}\emph{\cite{31}}\label{lem12}
Let $\mathcal{A}\in{\mathbb{R}}^{[m,n]}$ be a nonnegative tensor. Then $\rho(\mathcal{A})\geq\max\limits_{i\in{N}}\{a_{i\cdots i}\}$.
\end{lem12}
\newtheorem{lem4}[lem1]{Lemma}
\begin{lem4}\emph{\cite{6}}\label{lem3}
Let $\mathcal{A}$ be a strong M-tensor and denoted by $\tau(\mathcal{A})$ the minimum value of the real part of all eigenvalues of $\mathcal{A}$. Then $\tau(\mathcal{A})>0$ is an eigenvalue of $\mathcal{A}$ with a nonnegative eigenvector. Moreover, if $\mathcal{A}$ is irreducible, then $\tau(\mathcal{A})$ is a unique eigenvalue with a positive eigenvector.
\end{lem4}
\newtheorem{lem5}[lem1]{Lemma}
\begin{lem5}\emph{\cite{7}}\label{lem3}
Let $\mathcal{A}$ be a weakly irreducible strong M-tensor. Then $\tau(\mathcal{A})\leq\min\limits_{i\in{N}}\{a_{i\cdots i}\}$.
\end{lem5}
\newtheorem{lem7}[lem1]{Lemma}
\begin{lem7}\emph{\cite{16}}\label{lem3}
Let $a,b,c\geq{0}$ and $d>0$.\\
(I)~If $\frac{a}{b+c+d}\leq{1}$, then
\begin{eqnarray}
\frac{a-(b+c)}{d}\leq\frac{a-b}{c+d}\leq\frac{a}{b+c+d}.\nonumber
\end{eqnarray}
(II)~If $\frac{a}{b+c+d}\geq{1}$, then
\begin{eqnarray}
\frac{a-(b+c)}{d}\geq\frac{a-b}{c+d}\geq\frac{a}{b+c+d}.\nonumber
\end{eqnarray}
\end{lem7}

\section{A new \emph{S}-type eigenvalue localization set for tensors}
In this section, we investigate eigenvalue localization sets and present a new \emph{S}-type eigenvalue localization set for tensors, and the comparison results of this new set with those in Lemmas 1.1-1.4 are established.
\newtheorem{thm1}{Theorem}[section]
\begin{thm1}\label{thm1}
Let $\mathcal{A}=(a_{i_{1}\cdots i_{m}})\in{\mathbb{C}}^{[m,n]}$, $n\geq{2}$ and $S$ be a nonempty proper subset of $N$. Then
\begin{eqnarray}
\sigma(\mathcal{A})\subseteq\Upsilon^{S}(\mathcal{A}):=\Big(\Upsilon_{i,j}^{S}(\mathcal{A})\Big)\bigcup
\Big(\Upsilon_{i,j}^{\bar{S}}(\mathcal{A})\Big),
\end{eqnarray}
where
\begin{eqnarray}
&&\Upsilon_{i,j}^{S}(\mathcal{A})=\left(\bigcup_{i\in{S}}\hat{\Upsilon}_{i}^{1}(\mathcal{A})\right)\bigcup\left(\bigcup_{i\in{S},j\in{\bar{S}}}\Big(\tilde{\Upsilon}_{i,j}^{1}(\mathcal{A})\bigcap\Gamma_{i}(\mathcal{A})\Big)\right),\nonumber\\
&&\Upsilon_{i,j}^{\bar{S}}(\mathcal{A})=\left(\bigcup_{i\in{\bar{S}}}\hat{\Upsilon}_{i}^{2}(\mathcal{A})\right)\bigcup\left(\bigcup_{i\in{\bar{S}},j\in{S}}\Big(\tilde{\Upsilon}_{i,j}^{2}(\mathcal{A})\bigcap\Gamma_{i}(\mathcal{A})\Big)\right),\nonumber
\end{eqnarray}
with
\begin{eqnarray}
&&\hat{\Upsilon}_{i}^{1}(\mathcal{A})=\{z\in{\mathbb{C}}:|z-a_{i\cdots i}|\leq r_{i}^{\overline{\Delta^{\bar{S}}}}(\mathcal{A})\},\nonumber\\
&&\hat{\Upsilon}_{i}^{2}(\mathcal{A})=\{z\in{\mathbb{C}}:|z-a_{i\cdots i}|\leq r_{i}^{\overline{\Delta^{S}}}(\mathcal{A})\},\nonumber\\
&&\tilde{\Upsilon}_{i,j}^{1}(\mathcal{A})=\{z\in{\mathbb{C}}:(|z-a_{i\cdots i}|-r_{i}^{\overline{\Delta^{\bar{S}}}}(\mathcal{A}))(|z-a_{j\cdots j}|-r_{j}^{\Delta^{\bar{S}}}(\mathcal{A}))
\leq r_{i}^{\Delta^{\bar{S}}}(\mathcal{A})r_{j}^{\overline{\Delta^{\bar{S}}}}(\mathcal{A})\},\nonumber\\
&&\tilde{\Upsilon}_{i,j}^{2}(\mathcal{A})=\{z\in{\mathbb{C}}:(|z-a_{i\cdots i}|-r_{i}^{\overline{\Delta^{{S}}}}(\mathcal{A}))(|z-a_{j\cdots j}|-r_{j}^{\Delta^{S}}(\mathcal{A}))\leq r_{i}^{\Delta^{S}}(\mathcal{A})r_{j}^{\overline{\Delta^{{S}}}}(\mathcal{A})\}.\nonumber
\end{eqnarray}
\end{thm1}
\noindent {\bf Proof.} For any $\lambda\in\sigma(\mathcal{A})$, let $x=(x_{1},x_{2},\cdots,x_{n})^{T}\in{\mathbb{C}}^{n}\backslash{0}$ be an associated eigenvector, i.e.,
\begin{eqnarray}
\mathcal{A}x^{m-1}=\lambda x^{[m-1]}.
\end{eqnarray}
Let $|x_{p}|=\max\limits_{i\in{S}}\{|x_{i}|\}$ and $|x_{q}|=\max\limits_{i\in{\bar{S}}}\{|x_{i}|\}$. Then, $x_{p}\neq{0}$ or $x_{q}\neq{0}$.  Now, let us distinguish two cases to prove.

(i)~$|x_{p}|\geq|x_{q}|$, so $|x_{p}|=\max\limits_{i\in{N}}\{|x_{i}|\}$ and $|x_{p}|>0$. It follows from (2) that
\begin{eqnarray}
\sum\limits_{i_{2},\cdots, i_{m}=1}^{n}a_{pi_{2}\cdots i_{m}}x_{i_{2}}\cdots x_{i_{m}}=\lambda x_{p}^{m-1}.\nonumber
\end{eqnarray}
Hence, we have
\begin{eqnarray}
(\lambda-a_{p\cdots p})x_{p}^{m-1}=\sum\limits_{\substack{(i_{2},\cdots,i_{m})\in{\overline{\Delta^{\bar{S}}}},\\ \delta_{pi_{2}\cdots i_{m}=0}}}a_{pi_{2}\cdots i_{m}}x_{i_{2}}\cdots x_{i_{m}}+\sum\limits_{(i_{2},\cdots,i_{m})\in\Delta^{\bar{S}}}a_{pi_{2}\cdots i_{m}}x_{i_{2}}\cdots x_{i_{m}}.\nonumber
\end{eqnarray}
Taking absolute values in the above equation and using the triangle inequality yield
\begin{eqnarray}
&&|\lambda-a_{p\cdots p}||x_{p}|^{m-1}\nonumber\\&\leq&\sum\limits_{\substack{(i_{2},\cdots,i_{m})\in{\overline{\Delta^{\bar{S}}}},\\ \delta_{pi_{2}\cdots i_{m}=0}}}|a_{pi_{2}\cdots i_{m}}||x_{i_{2}}|\cdots |x_{i_{m}}|+\sum\limits_{(i_{2},\cdots,i_{m})\in\Delta^{\bar{S}}}|a_{pi_{2}\cdots i_{m}}||x_{i_{2}}|\cdots |x_{i_{m}}|\nonumber\\
&\leq&\sum\limits_{\substack{(i_{2},\cdots,i_{m})\in{\overline{\Delta^{\bar{S}}}},\\ \delta_{pi_{2}\cdots i_{m}=0}}}|a_{pi_{2}\cdots i_{m}}||x_{p}|^{m-1}+\sum\limits_{(i_{2},\cdots,i_{m})\in\Delta^{\bar{S}}}|a_{pi_{2}\cdots i_{m}}||x_{q}|^{m-1}\nonumber\\
&=&r_{p}^{\overline{\Delta^{\bar{S}}}}(\mathcal{A})|x_{p}|^{m-1}+r_{p}^{\Delta^{\bar{S}}}(\mathcal{A})|x_{q}|^{m-1},\nonumber
\end{eqnarray}
which means that
\begin{eqnarray}
(|\lambda-a_{p\cdots p}|-r_{p}^{\overline{\Delta^{\bar{S}}}}(\mathcal{A}))|x_{p}|^{m-1}\leq r_{p}^{\Delta^{\bar{S}}}(\mathcal{A})|x_{q}|^{m-1}.
\end{eqnarray}
If $|x_{q}|=0$, it follows from (3) that $|\lambda-a_{p\cdots p}|-r_{p}^{\overline{\Delta^{\bar{S}}}}(\mathcal{A})\leq 0$ by $|x_{p}|>0$, that is, $|\lambda-a_{p\cdots p}|\leq r_{p}^{\overline{\Delta^{\bar{S}}}}(\mathcal{A})$. Evidently, $\lambda\in\hat{\Upsilon}_{p}^{1}(\mathcal{A})\subseteq\Upsilon^{S}(\mathcal{A})$. Otherwise, $|x_{q}|>0$. If $\lambda\notin\bigcup\limits_{i\in{S}}\hat{\Upsilon}_{i}^{1}(\mathcal{A})$, it is easy to see that for any $i\in{S}$,
\begin{eqnarray}
|\lambda-a_{i\cdots i}|>r_{i}^{\overline{\Delta^{\bar{S}}}}(\mathcal{A}).\nonumber
\end{eqnarray}
In particular, $|\lambda-a_{p\cdots p}|>r_{p}^{\overline{\Delta^{\bar{S}}}}(\mathcal{A})$, i.e., $|\lambda-a_{p\cdots p}|-r_{p}^{\overline{\Delta^{\bar{S}}}}(\mathcal{A})>0$. By (3), it is not difficult to verify that $\lambda\in\Gamma_{p}(\mathcal{A})$. Besides, it follows from (2) that
\begin{eqnarray}
&&|\lambda-a_{q\cdots q}||x_{q}|^{m-1}\nonumber\\&\leq&\sum\limits_{(i_{2},\cdots,i_{m})\in{\overline{\Delta^{\bar{S}}}}}|a_{qi_{2}\cdots i_{m}}||x_{i_{2}}|\cdots |x_{i_{m}}|+\sum\limits_{\substack{(i_{2},\cdots,i_{m})\in{\Delta^{\bar{S}}},\\ \delta_{qi_{2}\cdots i_{m}=0}}}|a_{qi_{2}\cdots i_{m}}||x_{i_{2}}|\cdots |x_{i_{m}}|\nonumber\\
&\leq&\sum\limits_{(i_{2},\cdots,i_{m})\in{\overline{\Delta^{\bar{S}}}}}|a_{qi_{2}\cdots i_{m}}||x_{p}|^{m-1}+\sum\limits_{\substack{(i_{2},\cdots,i_{m})\in{\Delta^{\bar{S}}},\\ \delta_{qi_{2}\cdots i_{m}=0}}}|a_{qi_{2}\cdots i_{m}}||x_{q}|^{m-1}\nonumber\\
&=&r_{q}^{\overline{\Delta^{\bar{S}}}}(\mathcal{A})|x_{p}|^{m-1}+r_{q}^{\Delta^{\bar{S}}}(\mathcal{A})|x_{q}|^{m-1},\nonumber
\end{eqnarray}
which is equivalent to
\begin{eqnarray}
(|\lambda-a_{q\cdots q}|-r_{q}^{\Delta^{\bar{S}}}(\mathcal{A}))|x_{q}|^{m-1}\leq r_{q}^{\overline{\Delta^{\bar{S}}}}(\mathcal{A})|x_{p}|^{m-1}.
\end{eqnarray}
Note that $|x_{p}|>0$ and $|\lambda-a_{p\cdots p}|>r_{p}^{\overline{\Delta^{\bar{S}}}}(\mathcal{A})$, multiplying (3) with (4) results in
\begin{eqnarray}
&&(|\lambda-a_{p\cdots p}|-r_{p}^{\overline{\Delta^{\bar{S}}}}(\mathcal{A}))(|\lambda-a_{q\cdots q}|-r_{q}^{\Delta^{\bar{S}}}(\mathcal{A}))|x_{q}|^{m-1}|x_{p}|^{m-1}\nonumber\\
&\leq&r_{p}^{\Delta^{\bar{S}}}(\mathcal{A})r_{q}^{\overline{\Delta^{\bar{S}}}}(\mathcal{A})|x_{p}|^{m-1}|x_{q}|^{m-1},\nonumber
\end{eqnarray}
which implies that
\begin{eqnarray}
(|\lambda-a_{p\cdots p}|-r_{p}^{\overline{\Delta^{\bar{S}}}}(\mathcal{A}))(|\lambda-a_{q\cdots q}|-r_{q}^{\Delta^{\bar{S}}}(\mathcal{A}))
\leq r_{p}^{\Delta^{\bar{S}}}(\mathcal{A})r_{q}^{\overline{\Delta^{\bar{S}}}}(\mathcal{A})\nonumber
\end{eqnarray}
by $|x_{p}|\geq|x_{q}|>0$. Therefore, $\lambda\in\Big(\tilde{\Upsilon}_{p,q}^{1}(\mathcal{A})\bigcap\Gamma_{p}(\mathcal{A})\Big)\subseteq\Upsilon^{S}(\mathcal{A})$.

(ii)~$|x_{p}|\leq|x_{q}|$, so $|x_{q}|=\max\limits_{i\in{N}}\{|x_{i}|\}$ and $|x_{q}|>0$. It follows from (2) that
\begin{eqnarray}
\sum\limits_{i_{2},\cdots, i_{m}=1}^{n}a_{qi_{2}\cdots i_{m}}x_{i_{2}}\cdots x_{i_{m}}=\lambda x_{q}^{m-1}.\nonumber
\end{eqnarray}
Therefore, we have
\begin{eqnarray}
(\lambda-a_{q\cdots q})x_{q}^{m-1}=\sum\limits_{(i_{2},\cdots,i_{m})\in{\Delta^{{S}}}}a_{qi_{2}\cdots i_{m}}x_{i_{2}}\cdots x_{i_{m}}+\sum\limits_{\substack{(i_{2},\cdots,i_{m})\in\overline{\Delta^{S}},\\ \delta_{qi_{2}\cdots i_{m}=0}}}a_{qi_{2}\cdots i_{m}}x_{i_{2}}\cdots x_{i_{m}}.\nonumber
\end{eqnarray}
Taking modulus in the above equation and using the triangle inequality give
\begin{eqnarray}
&&|\lambda-a_{q\cdots q}||x_{q}|^{m-1}\nonumber\\&\leq&\sum\limits_{(i_{2},\cdots,i_{m})\in{\Delta^{{S}}}}|a_{qi_{2}\cdots i_{m}}||x_{i_{2}}|\cdots |x_{i_{m}}|+\sum\limits_{\substack{(i_{2},\cdots,i_{m})\in\overline{\Delta^{S}},\\ \delta_{qi_{2}\cdots i_{m}=0}}}|a_{qi_{2}\cdots i_{m}}||x_{i_{2}}|\cdots |x_{i_{m}}|\nonumber\\
&\leq&\sum\limits_{(i_{2},\cdots,i_{m})\in{\Delta^{{S}}}}|a_{qi_{2}\cdots i_{m}}||x_{p}|^{m-1}+\sum\limits_{\substack{(i_{2},\cdots,i_{m})\in\overline{\Delta^{S}},\\ \delta_{qi_{2}\cdots i_{m}=0}}}|a_{qi_{2}\cdots i_{m}}||x_{q}|^{m-1}\nonumber\\
&=&r_{q}^{\Delta^{{S}}}(\mathcal{A})|x_{p}|^{m-1}+r_{q}^{\overline{\Delta^{S}}}(\mathcal{A})|x_{q}|^{m-1},\nonumber
\end{eqnarray}
which yields that
\begin{eqnarray}
(|\lambda-a_{q\cdots q}|-r_{q}^{\overline{\Delta^{S}}}(\mathcal{A}))|x_{q}|^{m-1}\leq r_{q}^{\Delta^{{S}}}(\mathcal{A})|x_{p}|^{m-1}.
\end{eqnarray}
If $|x_{p}|=0$, it follows from (5) that $|\lambda-a_{q\cdots q}|-r_{q}^{\overline{\Delta^{S}}}(\mathcal{A})\leq 0$ by $|x_{q}|>0$, i.e., $|\lambda-a_{q\cdots q}|\leq r_{q}^{\overline{\Delta^{S}}}(\mathcal{A})$, obviously, $\lambda\in\hat{\Upsilon}_{q}^{2}(\mathcal{A})\subseteq\Upsilon^{S}(\mathcal{A})$. Otherwise, $|x_{p}|>0$. If $\lambda\notin\bigcup\limits_{i\in{\bar{S}}}\hat{\Upsilon}_{i}^{2}(\mathcal{A})$, we are easy to see that for any $i\in{\bar{S}}$,
\begin{eqnarray}
|\lambda-a_{i\cdots i}|>r_{i}^{\overline{\Delta^{S}}}(\mathcal{A}).\nonumber
\end{eqnarray}
In particular, $|\lambda-a_{q\cdots q}|>r_{q}^{\overline{\Delta^{S}}}(\mathcal{A})$, i.e., $|\lambda-a_{q\cdots q}|-r_{q}^{\overline{\Delta^{S}}}(\mathcal{A})>0$. By (5), we infer that $\lambda\in\Gamma_{q}(\mathcal{A})$. In addition, it follows from (2) that
\begin{eqnarray}
&&|\lambda-a_{p\cdots p}||x_{p}|^{m-1}\nonumber\\&\leq&\sum\limits_{\substack{(i_{2},\cdots,i_{m})\in{\Delta^{S}},\\ \delta_{pi_{2}\cdots i_{m}=0}}}|a_{pi_{2}\cdots i_{m}}||x_{i_{2}}|\cdots |x_{i_{m}}|+\sum\limits_{(i_{2},\cdots,i_{m})\in{\overline{\Delta^{S}}}}|a_{pi_{2}\cdots i_{m}}||x_{i_{2}}|\cdots |x_{i_{m}}|\nonumber\\
&\leq&\sum\limits_{\substack{(i_{2},\cdots,i_{m})\in{\Delta^{S}},\\ \delta_{pi_{2}\cdots i_{m}=0}}}|a_{pi_{2}\cdots i_{m}}||x_{p}|^{m-1}+\sum\limits_{(i_{2},\cdots,i_{m})\in{\overline{\Delta^{S}}}}|a_{pi_{2}\cdots i_{m}}||x_{q}|^{m-1}\nonumber\\
&=&r_{p}^{\Delta^{S}}(\mathcal{A})|x_{p}|^{m-1}+r_{p}^{\overline{\Delta^{S}}}(\mathcal{A})|x_{q}|^{m-1},\nonumber
\end{eqnarray}
which is equivalent to
\begin{eqnarray}
(|\lambda-a_{p\cdots p}|-r_{p}^{\Delta^{S}}(\mathcal{A}))|x_{p}|^{m-1}\leq r_{p}^{\overline{\Delta^{S}}}(\mathcal{A})|x_{q}|^{m-1}.
\end{eqnarray}
Having in mind that $|x_{q}|>0$ and $|\lambda-a_{q\cdots q}|>r_{q}^{\overline{\Delta^{S}}}(\mathcal{A})$, multiplying (5) with (6) results in
\begin{eqnarray}
&&(|\lambda-a_{p\cdots p}|-r_{p}^{\Delta^{S}}(\mathcal{A}))(|\lambda-a_{q\cdots q}|-r_{q}^{\overline{\Delta^{S}}}(\mathcal{A}))|x_{q}|^{m-1}|x_{p}|^{m-1}\nonumber\\
&\leq&r_{p}^{\overline{\Delta^{S}}}(\mathcal{A})r_{q}^{{\Delta^{S}}}(\mathcal{A})|x_{p}|^{m-1}|x_{q}|^{m-1},\nonumber
\end{eqnarray}
which results in
\begin{eqnarray}
(|\lambda-a_{p\cdots p}|-r_{p}^{\Delta^{S}}(\mathcal{A}))(|\lambda-a_{q\cdots q}|-r_{q}^{\overline{\Delta^{S}}}(\mathcal{A}))
\leq r_{p}^{\overline{\Delta^{S}}}(\mathcal{A})r_{q}^{{\Delta^{S}}}(\mathcal{A})\nonumber
\end{eqnarray}
by $|x_{q}|\geq|x_{p}|>0$. This leads to $\lambda\in\Big(\tilde{\Upsilon}_{q,p}^{2}(\mathcal{A})\bigcap\Gamma_{q}(\mathcal{A})\Big)\subseteq\Upsilon^{S}(\mathcal{A})$.
This completes our proof of Theorem 3.1. \hfill$\blacksquare$

Now, we establish a comparison result between $\Upsilon^{S}(\mathcal{A})$, $\Omega^{S}(\mathcal{A})$, $\mathcal{K}^{S}(\mathcal{A})$, $\mathcal{K}(\mathcal{A})$ and $\Gamma(\mathcal{A})$ as follows.
\newtheorem{thm11}[thm1]{Theorem}
\begin{thm11}\label{thm11}
Let $\mathcal{A}=(a_{i_{1}\cdots i_{m}})\in{\mathbb{C}}^{[m,n]}$, $n\geq{2}$ and $S$ be a nonempty proper subset of $N$. Then
\begin{eqnarray}
\Upsilon^{S}(\mathcal{A})\subseteq\Omega^{S}(\mathcal{A})\subseteq\mathcal{K}^{S}(\mathcal{A})\subseteq\mathcal{K}(\mathcal{A})\subseteq\Gamma(\mathcal{A}).\nonumber
\end{eqnarray}
\end{thm11}
\noindent {\bf Proof.} By Theorem 6 in \cite{16}, we see that $\Omega^{S}(\mathcal{A})\subseteq\mathcal{K}^{S}(\mathcal{A})\subseteq\mathcal{K}(\mathcal{A})\subseteq\Gamma(\mathcal{A})$ holds. Thus, we only need to prove $\Upsilon^{S}(\mathcal{A})\subseteq\Omega^{S}(\mathcal{A})$. Let $z\in\Upsilon^{S}(\mathcal{A})$. Then
\begin{eqnarray}
z\in\Upsilon_{i,j}^{S}(\mathcal{A})\ \mathrm{or}\ z\in\Upsilon_{i,j}^{\bar{S}}(\mathcal{A}).\nonumber
\end{eqnarray}
Without loss of generality, we first assume that $z\in\Upsilon_{i,j}^{S}(\mathcal{A})$.
If $z\in\bigcup\limits_{i\in{S}}\hat{\Upsilon}_{i}^{1}(\mathcal{A})$, then there exists one index $i_{0}\in{S}$ such that
\begin{eqnarray}
|z-a_{i_{0}\cdots i_{0}}|\leq r_{i_{0}}^{\overline{\Delta^{\bar{S}}}}(\mathcal{A}),\nonumber
\end{eqnarray}
i.e., $|z-a_{i_{0}\cdots i_{0}}|-r_{i_{0}}^{\overline{\Delta^{\bar{S}}}}(\mathcal{A})\leq{0}$. Hence, for any $i\in{\bar{S}}$, it follows
\begin{eqnarray}
(|z-a_{i\cdots i}|)(|z-a_{i_{0}\cdots i_{0}}|-r_{i_{0}}^{\overline{\Delta^{\bar{S}}}}(\mathcal{A}))\leq r_{i}(\mathcal{A})r_{i_{0}}^{\Delta^{\bar{S}}}(\mathcal{A}),\nonumber
\end{eqnarray}
which implies that $z\in\Omega_{i,i_{0}}^{\bar{S}}(\mathcal{A})\subseteq\Omega^{S}(\mathcal{A})$. Otherwise, $z\notin\bigcup\limits_{i\in{S}}\hat{\Upsilon}_{i}^{1}(\mathcal{A})$, then
\begin{eqnarray}
z\in\left(\bigcup_{i\in{S},j\in{\bar{S}}}\Big(\tilde{\Upsilon}_{i,j}^{1}(\mathcal{A})\bigcap\Gamma_{i}(\mathcal{A})\Big)\right)
\end{eqnarray}
and
\begin{eqnarray}
|z-a_{i\cdots i}|>r_{i}^{\overline{\Delta^{\bar{S}}}}(\mathcal{A})
\end{eqnarray}
for any $i\in{S}$. It follows from (7) that there exist $p\in{S}$ and $q\in{\bar{S}}$ such that
\begin{eqnarray}
|z-a_{p\cdots p}|\leq r_{p}(\mathcal{A})
\end{eqnarray}
and
\begin{eqnarray}
(|z-a_{p\cdots p}|-r_{p}^{\overline{\Delta^{\bar{S}}}}(\mathcal{A}))(|z-a_{q\cdots q}|-r_{q}^{\Delta^{\bar{S}}}(\mathcal{A}))
\leq r_{p}^{\Delta^{\bar{S}}}(\mathcal{A})r_{q}^{\overline{\Delta^{\bar{S}}}}(\mathcal{A}).
\end{eqnarray}
If $r_{p}^{\Delta^{\bar{S}}}(\mathcal{A})r_{q}^{\overline{\Delta^{\bar{S}}}}(\mathcal{A})=0$, combining (8) and (10) results in
\begin{eqnarray}
|z-a_{q\cdots q}|-r_{q}^{\Delta^{\bar{S}}}(\mathcal{A})\leq{0}\leq r_{q}^{\overline{\Delta^{\bar{S}}}}(\mathcal{A}),
\end{eqnarray}
that is, $|z-a_{q\cdots q}|\leq{r_{q}(\mathcal{A})}$, which is equivalent to
\begin{eqnarray}
|z-a_{q\cdots q}|-r_{q}^{\overline{\Delta^{{S}}}}(\mathcal{A})\leq r_{q}^{\Delta^{S}}(\mathcal{A}).
\end{eqnarray}
Multiplying (9) with (12) yields
\begin{eqnarray}
(|z-a_{p\cdots p}|)(|z-a_{q\cdots q}|-r_{q}^{\overline{\Delta^{S}}}(\mathcal{A}))\leq r_{p}(\mathcal{A})r_{q}^{\Delta^{S}}(\mathcal{A}).
\end{eqnarray}
This means that $z\in\Omega_{p,q}^{S}(\mathcal{A})\subseteq\Omega^{S}(\mathcal{A})$.

In the sequel, we discuss the case $r_{p}^{\Delta^{\bar{S}}}(\mathcal{A})r_{q}^{\overline{\Delta^{\bar{S}}}}(\mathcal{A})>0$, then by dividing (10) by $r_{p}^{\Delta^{\bar{S}}}(\mathcal{A})r_{q}^{\overline{\Delta^{\bar{S}}}}(\mathcal{A})$ is given by
\begin{eqnarray}
\frac{(|z-a_{p\cdots p}|-r_{p}^{\overline{\Delta^{\bar{S}}}}(\mathcal{A}))(|z-a_{q\cdots q}|-r_{q}^{\Delta^{\bar{S}}}(\mathcal{A}))}{r_{p}^{\Delta^{\bar{S}}}(\mathcal{A})r_{q}^{\overline{\Delta^{\bar{S}}}}(\mathcal{A})}\leq 1.
\end{eqnarray}
If $\frac{|z-a_{q\cdots q}|-r_{q}^{\Delta^{\bar{S}}}(\mathcal{A})}{r_{q}^{\overline{\Delta^{\bar{S}}}}(\mathcal{A})}\geq{1}$, let $a=|z-a_{q\cdots q}|\geq{0}$, $b+c=r_{q}^{\Delta^{\bar{S}}}(\mathcal{A})\geq{0}$ with $b,c\geq{0}$ and $d=r_{q}^{\overline{\Delta^{\bar{S}}}}(\mathcal{A})>0$, then by (II) of Lemma 2.5, we have
\begin{eqnarray}
\frac{|z-a_{p\cdots p}|-r_{p}^{\overline{\Delta^{\bar{S}}}}(\mathcal{A})}{r_{p}^{{\Delta^{\bar{S}}}}(\mathcal{A})}\frac{|z-a_{q\cdots q}|}{r_{q}(\mathcal{A})}\leq
\frac{|z-a_{p\cdots p}|-r_{p}^{\overline{\Delta^{\bar{S}}}}(\mathcal{A})}{r_{p}^{\Delta^{\bar{S}}}(\mathcal{A})}\frac{|z-a_{q\cdots q}|-r_{q}^{\Delta^{\bar{S}}}(\mathcal{A})}{r_{q}^{\overline{\Delta^{\bar{S}}}}(\mathcal{A})}\leq 1,\nonumber
\end{eqnarray}
which is equivalent to
\begin{eqnarray}
(|z-a_{q\cdots q}|)(|z-a_{p\cdots p}|-r_{p}^{\overline{\Delta^{\bar{S}}}}(\mathcal{A}))\leq r_{q}(\mathcal{A})r_{p}^{\Delta^{\bar{S}}}(\mathcal{A}).
\end{eqnarray}
This implies that $z\in\Omega_{q,p}^{\bar{S}}(\mathcal{A})\subseteq\Omega^{S}(\mathcal{A})$. Furthermore, if $\frac{|z-a_{q\cdots q}|-r_{q}^{\Delta^{\bar{S}}}(\mathcal{A})}{r_{q}^{\overline{\Delta^{\bar{S}}}}(\mathcal{A})}\leq{1}$, then
(12) holds. Multiplying (9) with (12) leads to
\begin{eqnarray}
(|z-a_{p\cdots p}|)(|z-a_{q\cdots q}|-r_{q}^{\overline{\Delta^{S}}}(\mathcal{A}))\leq r_{p}(\mathcal{A})r_{q}^{\Delta^{S}}(\mathcal{A}),
\end{eqnarray}
which implies that $z\in\Omega_{p,q}^{S}(\mathcal{A})\subseteq\Omega^{S}(\mathcal{A})$.

On the other hand, we prove the case $z\in\Upsilon_{i,j}^{\bar{S}}(\mathcal{A})$. If $z\in\bigcup\limits_{i\in{\bar{S}}}\hat{\Upsilon}_{i}^{2}(\mathcal{A})$, then there is one index $i_{1}\in{\bar{S}}$ such that
\begin{eqnarray}
|z-a_{i_{1}\cdots i_{1}}|\leq r_{i_{1}}^{\overline{\Delta^{S}}}(\mathcal{A}),\nonumber
\end{eqnarray}
i.e., $|z-a_{i_{1}\cdots i_{1}}|-r_{i_{1}}^{\overline{\Delta^{S}}}(\mathcal{A})\leq{0}$. Then, for any $i\in{S}$, we deduce that
\begin{eqnarray}
(|z-a_{i\cdots i}|)(|z-a_{i_{1}\cdots i_{1}}|-r_{i_{1}}^{\overline{\Delta^{S}}}(\mathcal{A}))\leq r_{i}(\mathcal{A})r_{i_{1}}^{\Delta^{S}}(\mathcal{A}),\nonumber
\end{eqnarray}
which means that $z\in\Omega_{i,i_{1}}^{S}(\mathcal{A})\subseteq\Omega^{S}(\mathcal{A})$. In addition, $z\notin\bigcup\limits_{i\in{\bar{S}}}\hat{\Upsilon}_{i}^{2}(\mathcal{A})$, then
\begin{eqnarray}
z\in\left(\bigcup_{i\in{\bar{S}},j\in{S}}\Big(\tilde{\Upsilon}_{i,j}^{2}(\mathcal{A})\bigcap\Gamma_{i}(\mathcal{A})\Big)\right)
\end{eqnarray}
and
\begin{eqnarray}
|z-a_{j\cdots j}|>r_{j}^{\overline{\Delta^{S}}}(\mathcal{A})
\end{eqnarray}
for any $j\in{\bar{S}}$. It follows from (17) that there exist $p\in{S}$ and $q\in{\bar{S}}$ such that
\begin{eqnarray}
|z-a_{q\cdots q}|\leq r_{q}(\mathcal{A})
\end{eqnarray}
and
\begin{eqnarray}
(|z-a_{q\cdots q}|-r_{q}^{\overline{\Delta^{{S}}}}(\mathcal{A}))(|z-a_{p\cdots p}|-r_{p}^{\Delta^{S}}(\mathcal{A}))\leq r_{q}^{\Delta^{S}}(\mathcal{A})r_{p}^{\overline{\Delta^{{S}}}}(\mathcal{A}).
\end{eqnarray}
If $r_{q}^{\Delta^{S}}(\mathcal{A})r_{p}^{\overline{\Delta^{{S}}}}(\mathcal{A})=0$, combining (18) and (20) results in
\begin{eqnarray}
|z-a_{p\cdots p}|-r_{p}^{\Delta^{S}}(\mathcal{A})\leq{0}\leq r_{p}^{\overline{\Delta^{S}}}(\mathcal{A}),
\end{eqnarray}
which leads to $|z-a_{p\cdots p}|\leq{r_{p}(\mathcal{A})}$, and therefore
\begin{eqnarray}
|z-a_{p\cdots p}|-r_{p}^{\overline{\Delta^{{\bar{S}}}}}(\mathcal{A})\leq r_{p}^{\Delta^{\bar{S}}}(\mathcal{A}).
\end{eqnarray}
Multiplying (19) with (22) derives
\begin{eqnarray}
(|z-a_{q\cdots q}|)(|z-a_{p\cdots p}|-r_{p}^{\overline{\Delta^{\bar{S}}}}(\mathcal{A}))\leq r_{q}(\mathcal{A})r_{p}^{\Delta^{\bar{S}}}(\mathcal{A}).
\end{eqnarray}
It follows from (23) that $z\in\Omega_{q,p}^{\bar{S}}(\mathcal{A})\subseteq\Omega^{S}(\mathcal{A})$.

Afterwards, we investigate the case $r_{q}^{\Delta^{S}}(\mathcal{A})r_{p}^{\overline{\Delta^{{S}}}}(\mathcal{A})>0$, then by dividing (20) by $r_{q}^{\Delta^{S}}(\mathcal{A})r_{p}^{\overline{\Delta^{{S}}}}(\mathcal{A})$, we have
\begin{eqnarray}
\frac{(|z-a_{q\cdots q}|-r_{q}^{\overline{\Delta^{{S}}}}(\mathcal{A}))(|z-a_{p\cdots p}|-r_{p}^{\Delta^{S}}(\mathcal{A}))}{r_{q}^{\Delta^{S}}(\mathcal{A})r_{p}^{\overline{\Delta^{{S}}}}(\mathcal{A})}\leq 1.
\end{eqnarray}
If $\frac{|z-a_{p\cdots p}|-r_{p}^{\Delta^{S}}(\mathcal{A})}{r_{p}^{\overline{\Delta^{{S}}}}(\mathcal{A})}\geq{1}$, let $a=|z-a_{p\cdots p}|\geq{0}$, $b+c=r_{p}^{\Delta^{S}}(\mathcal{A})\geq{0}$ with $b,c\geq{0}$ and $d=r_{p}^{\overline{\Delta^{{S}}}}(\mathcal{A})>0$, then by (II) of Lemma 2.5, we have
\begin{eqnarray}
\frac{|z-a_{q\cdots q}|-r_{q}^{\overline{\Delta^{{S}}}}(\mathcal{A})}{r_{q}^{\Delta^{S}}(\mathcal{A})}\frac{|z-a_{p\cdots p}|}{r_{p}(\mathcal{A})}\leq
\frac{|z-a_{q\cdots q}|-r_{q}^{\overline{\Delta^{{S}}}}(\mathcal{A})}{r_{q}^{\Delta^{S}}(\mathcal{A})}\frac{|z-a_{p\cdots p}|-r_{p}^{\Delta^{S}}(\mathcal{A})}{r_{p}^{\overline{\Delta^{{S}}}}(\mathcal{A})}\leq 1,\nonumber
\end{eqnarray}
equivalently,
\begin{eqnarray}
(|z-a_{p\cdots p}|)(|z-a_{q\cdots q}|-r_{q}^{\overline{\Delta^{{S}}}}(\mathcal{A}))\leq r_{p}(\mathcal{A})r_{q}^{\Delta^{S}}(\mathcal{A}).
\end{eqnarray}
This implies that $z\in\Omega_{p,q}^{S}(\mathcal{A})\subseteq\Omega^{S}(\mathcal{A})$. Furthermore, if $\frac{|z-a_{p\cdots p}|-r_{p}^{\Delta^{S}}(\mathcal{A})}{r_{p}^{\overline{\Delta^{{S}}}}(\mathcal{A})}\leq{1}$, then
(22) holds. Multiplying (19) with (22) leads to
\begin{eqnarray}
(|z-a_{q\cdots q}|)(|z-a_{p\cdots p}|-r_{q}^{\overline{\Delta^{\bar{S}}}}(\mathcal{A}))\leq r_{q}(\mathcal{A})r_{p}^{\Delta^{\bar{S}}}(\mathcal{A}),
\end{eqnarray}
which implies that $z\in\Omega_{q,p}^{\bar{S}}(\mathcal{A})\subseteq\Omega^{S}(\mathcal{A})$.

It follows from the above discussions that $\Upsilon^{S}(\mathcal{A})\subseteq\Omega^{S}(\mathcal{A})$. The conclusion follows immediately from what we have proved. \hfill$\blacksquare$

\newtheorem{rem1}{Remark}[section]
\begin{rem1}\label{rem1}
\emph{For a complex tensor $\mathcal{A}\in{\mathbb{C}}^{[m,n]}$, $n\geq2$, the set $\Omega^{S}(\mathcal{A})$ consists of $|S|(n-|S|)$ sets $\Omega_{i,j}^{S}(\mathcal{A})$ and $|S|(n-|S|)$ sets $\Omega_{i,j}^{\bar{S}}(\mathcal{A})$, where $S$ is a nonempty proper subset of $N$, and therefore $\Omega^{S}(\mathcal{A})$ contains $2|S|(n-|S|)$ sets. In addition, the set $\Upsilon^{S}(\mathcal{A})$ consists of $|S|(n-|S|)$ sets $\tilde{\Upsilon}_{i,j}^{1}(\mathcal{A})$, $|S|(n-|S|)$ sets $\tilde{\Upsilon}_{i,j}^{2}(\mathcal{A})$, $|S|$ sets $\hat{\Upsilon}_{i}^{1}(\mathcal{A})$, $n-|S|$ sets $\hat{\Upsilon}_{i}^{2}(\mathcal{A})$ and $n$ sets $\Gamma_{i}(\mathcal{A})$, then there are $2|S|(n-|S|)+2n$ sets contained in $\Upsilon^{S}(\mathcal{A})$.
Hence there are more computations to determine $\Upsilon^{S}(\mathcal{A})$ than $\Omega^{S}(\mathcal{A})$, while $\Upsilon^{S}(\mathcal{A})$ can capture all eigenvalues of $\mathcal{A}$ more precisely than $\Omega^{S}(\mathcal{A})$ as showed in Theorem 3.2.}

\emph{Based on the above discussions, how to choose $S$ to make $\Upsilon^{S}(\mathcal{A})$ as sharp as possible is very interesting and important. However, this work is difficult especially the dimension of the tensor $\mathcal{A}$ is large. At present, it is very difficult for us to research this problem, we will continue to study this problem in the future.}
\end{rem1}
\section{Sufficient conditions for positive (semi-)definiteness of tensors}
As applications of the results in Section 3, we provide some checkable sufficient conditions for the positive definiteness and positive semi-definiteness of tensors, respectively in this section. Furthermore, a numerical example is implemented to illustrate the superiority of these conditions to those derived in \cite{16,21,31}.
\newtheorem{thm3}[thm1]{Theorem}
\begin{thm3}\label{thm3}
Let $\mathcal{A}=(a_{i_{1}\cdots i_{m}})\in{\mathbb{R}}^{[m,n]}$ be an even-order symmetric tensor with $a_{k\cdots k}>0$ for all $k\in{N}$. If there is a nonempty proper subset $S$ of $N$ and the following four statements hold:\\
(i)~$a_{i\cdots i}>r_{i}^{\overline{\Delta^{\bar{S}}}}(\mathcal{A})$ for any $i\in{S}$;\\
(ii)~$a_{i\cdots i}>r_{i}^{\overline{\Delta^{S}}}(\mathcal{A})$ for any $i\in{\bar{S}}$;\\
(iii)~For any $i\in{S}$, $j\in{\bar{S}}$,
\begin{eqnarray}
(a_{i\cdots i}-r_{i}^{\overline{\Delta^{\bar{S}}}}(\mathcal{A}))(a_{j\cdots j}-r_{j}^{\Delta^{\bar{S}}}(\mathcal{A}))
> r_{i}^{\Delta^{\bar{S}}}(\mathcal{A})r_{j}^{\overline{\Delta^{\bar{S}}}}(\mathcal{A})\nonumber
\end{eqnarray}
or $a_{i\cdots i}>r_{i}(\mathcal{A})$;\\
(iv)~For any $i\in{\bar{S}}$, $j\in{S}$,
\begin{eqnarray}
(a_{i\cdots i}-r_{i}^{\overline{\Delta^{{S}}}}(\mathcal{A}))(a_{j\cdots j}-r_{j}^{\Delta^{S}}(\mathcal{A}))>r_{i}^{\Delta^{S}}(\mathcal{A})r_{j}^{\overline{\Delta^{{S}}}}(\mathcal{A})\nonumber
\end{eqnarray}
or $a_{i\cdots i}>r_{i}(\mathcal{A})$, then $\mathcal{A}$ is positive definite.
\end{thm3}
\noindent {\bf Proof.}
Let $\lambda$ be an \emph{H}-eigenvalue of $\mathcal{A}$. We prove this theorem by assuming that $\lambda\leq0$ and leading a contradiction. From Theorem 3.1, we have $\lambda\in\Upsilon^{S}(\mathcal{A})$, which implies that there are $i_{0},i_{1},i_{2}\in{S}$ and $j_{0},j_{1},j_{2}\in{\bar{S}}$ such that
\begin{eqnarray}
|\lambda-a_{i_{0}\cdots i_{0}}|\leq r_{i_{0}}^{\overline{\Delta^{\bar{S}}}}(\mathcal{A})\ \mathrm{or} \ |\lambda-a_{j_{0}\cdots j_{0}}|\leq r_{j_{0}}^{\overline{\Delta^{S}}}(\mathcal{A})\nonumber
\end{eqnarray}
or
\begin{eqnarray}
&&(|\lambda-a_{i_{1}\cdots i_{1}}|-r_{i_{1}}^{\overline{\Delta^{\bar{S}}}}(\mathcal{A}))(|\lambda-a_{j_{1}\cdots j_{1}}|-r_{j_{1}}^{\Delta^{\bar{S}}}(\mathcal{A}))\leq r_{i_{1}}^{\Delta^{\bar{S}}}(\mathcal{A})r_{j_{1}}^{\overline{\Delta^{\bar{S}}}}(\mathcal{A}),\nonumber\\
&&|\lambda-a_{i_{1}\cdots i_{1}}|\leq r_{i_{1}}(\mathcal{A})\nonumber
\end{eqnarray}
or
\begin{eqnarray}
&&(|\lambda-a_{j_{2}\cdots j_{2}}|-r_{j_{2}}^{\overline{\Delta^{S}}}(\mathcal{A}))(|\lambda-a_{i_{2}\cdots i_{2}}|-r_{i_{2}}^{\Delta^{S}}(\mathcal{A}))\leq r_{j_{2}}^{\Delta^{S}}(\mathcal{A})r_{i_{2}}^{\overline{\Delta^{S}}}(\mathcal{A}),\nonumber\\
&&|\lambda-a_{j_{2}\cdots j_{2}}|\leq r_{j_{2}}(\mathcal{A}).\nonumber
\end{eqnarray}
It follows from $a_{k\cdots k}>0$ for all $k\in{N}$ that
\begin{eqnarray}
|\lambda-a_{i_{0}\cdots i_{0}}|\geq a_{i_{0}\cdots i_{0}}>r_{i_{0}}^{\overline{\Delta^{\bar{S}}}}(\mathcal{A})\ \mathrm{and} \ |\lambda-a_{j_{0}\cdots j_{0}}|\geq a_{j_{0}\cdots j_{0}}>r_{j_{0}}^{\overline{\Delta^{S}}}(\mathcal{A})\nonumber
\end{eqnarray}
and
\setlength\arraycolsep{2pt}
\begin{eqnarray}
&&(|\lambda-a_{i_{1}\cdots i_{1}}|-r_{i_{1}}^{\overline{\Delta^{\bar{S}}}}(\mathcal{A}))(|\lambda-a_{j_{1}\cdots j_{1}}|-r_{j_{1}}^{\Delta^{\bar{S}}}(\mathcal{A}))\nonumber\\
&\geq& (a_{i_{1}\cdots i_{1}}-r_{i_{1}}^{\overline{\Delta^{\bar{S}}}}(\mathcal{A}))(a_{j_{1}\cdots j_{1}}-r_{j_{1}}^{\Delta^{\bar{S}}}(\mathcal{A}))>r_{i_{1}}^{\Delta^{\bar{S}}}(\mathcal{A})r_{j_{1}}^{\overline{\Delta^{\bar{S}}}}(\mathcal{A})\nonumber
\end{eqnarray}
or $|\lambda-a_{i_{1}\cdots i_{1}}|\geq a_{i_{1}\cdots i_{1}}>r_{i_{1}}(\mathcal{A})$;
and
\setlength\arraycolsep{2pt}
\begin{eqnarray}
&&(|\lambda-a_{j_{2}\cdots j_{2}}|-r_{j_{2}}^{\overline{\Delta^{{S}}}}(\mathcal{A}))(|\lambda-a_{i_{2}\cdots i_{2}}|-r_{i_{2}}^{\Delta^{S}}(\mathcal{A}))\nonumber\\
&\geq& (a_{j_{2}\cdots j_{2}}-r_{j_{2}}^{\overline{\Delta^{{S}}}}(\mathcal{A}))(a_{i_{2}\cdots i_{2}}-r_{i_{2}}^{\Delta^{S}}(\mathcal{A}))> r_{j_{2}}^{\Delta^{S}}(\mathcal{A})r_{i_{2}}^{\overline{\Delta^{{S}}}}(\mathcal{A})\nonumber
\end{eqnarray}
or $|\lambda-a_{j_{2}\cdots j_{2}}|\geq a_{j_{2}\cdots j_{2}}>r_{j_{2}}(\mathcal{A})$. These lead to a contradiction. Hence, $\lambda>0$, and $\mathcal{A}$ is positive definite. \hfill$\blacksquare$

With the similar manner applied in the proof of Theorem 4.1, we can prove that $\mathcal{A}$ is positive semi-definite in the following theorem.
\newtheorem{thm10}[thm1]{Theorem}
\begin{thm10}\label{thm3}
Let $\mathcal{A}=(a_{i_{1}\cdots i_{m}})\in{\mathbb{R}}^{[m,n]}$ be an even-order symmetric tensor with $a_{k\cdots k}\geq0$ for all $k\in{N}$. If there is a nonempty proper subset $S$ of $N$ and the following four statements hold:\\
(i)~$a_{i\cdots i}\geq r_{i}^{\overline{\Delta^{\bar{S}}}}(\mathcal{A})$ for any $i\in{S}$;\\
(ii)~$a_{i\cdots i}\geq r_{i}^{\overline{\Delta^{S}}}(\mathcal{A})$ for any $i\in{\bar{S}}$;\\
(iii)~For any $i\in{S}$, $j\in{\bar{S}}$,
\begin{eqnarray}
(a_{i\cdots i}-r_{i}^{\overline{\Delta^{\bar{S}}}}(\mathcal{A}))(a_{j\cdots j}-r_{j}^{\Delta^{\bar{S}}}(\mathcal{A}))
\geq r_{i}^{\Delta^{\bar{S}}}(\mathcal{A})r_{j}^{\overline{\Delta^{\bar{S}}}}(\mathcal{A})\nonumber
\end{eqnarray}
or $a_{i\cdots i}\geq r_{i}(\mathcal{A})$;\\
(iv)~For any $i\in{\bar{S}}$, $j\in{S}$,
\begin{eqnarray}
(a_{i\cdots i}-r_{i}^{\overline{\Delta^{{S}}}}(\mathcal{A}))(a_{j\cdots j}-r_{j}^{\Delta^{S}}(\mathcal{A}))\geq r_{i}^{\Delta^{S}}(\mathcal{A})r_{j}^{\overline{\Delta^{{S}}}}(\mathcal{A})\nonumber
\end{eqnarray}
or $a_{i\cdots i}\geq r_{i}(\mathcal{A})$, then $\mathcal{A}$ is positive semi-definite.
\end{thm10}

The advantages of the results of Theorem 4.1 will be stressed by the following numerical example.
\begin{Exa}\label{E2}
\setlength\arraycolsep{4pt}
\emph{Let $\mathcal{A}=(a_{ijkl})\in{\mathbb{R}}^{[4,3]}$ be a real symmetric tensor with elements defined as follows:}
\begin{eqnarray}
&&a_{1111}=5.2,\ a_{2222}=6,\ a_{3333}=3.3,\ a_{1112}=-0.1,\ a_{1113}=0.1,\nonumber\\
&&a_{1122}=-0.2,\ a_{1123}=-0.2,\ a_{1133}=0,\ a_{1222}=-0.1,\ a_{1223}=0.3,\nonumber\\
&&a_{1233}=0.1,\ a_{1333}=-0.2,\ a_{2223}=0.1,\ a_{2233}=-0.1,\ a_{2333}=0.2.\nonumber
\end{eqnarray}
\end{Exa}
After some calculations, we conclude that the tensor $\mathcal{A}$ can not meet the conditions of Theorem 3.2 in \cite{21}, and for any nonempty proper subset $S$ of $N$, Theorem 4.2 of \cite{31} and Theorem 7 of \cite{16} can not be applied to determine the positive definiteness of $\mathcal{A}$, while we choose the nonempty proper subset $S$ of $N$ is $S=\{1,2\}$, then $\bar{S}=\{3\}$, thus following results are easy to obtain
\begin{eqnarray}
&&a_{1111}=5.2>3.7=r_{1}^{\overline{\Delta^{\bar{S}}}}(\mathcal{A}),\ a_{2222}=6>4.3=r_{2}^{\overline{\Delta^{\bar{S}}}}(\mathcal{A}),\nonumber\\
&&a_{3333}=3.3>2.1=r_{3}^{\overline{\Delta^{S}}}(\mathcal{A}),\nonumber\\
&&a_{1111}=5.2>3.9=r_{1}(\mathcal{A}),\ a_{2222}=6>4.5=r_{2}(\mathcal{A}),\nonumber\\
&&(a_{3333}-r_{3}^{\overline{\Delta^{{S}}}}(\mathcal{A}))(a_{1111}-r_{1}^{\Delta^{S}}(\mathcal{A}))=5.04>4.93=r_{3}^{\Delta^{S}}(\mathcal{A})r_{1}^{\overline{\Delta^{{S}}}}(\mathcal{A}),\nonumber\\
&&(a_{3333}-r_{3}^{\overline{\Delta^{{S}}}}(\mathcal{A}))(a_{2222}-r_{2}^{\Delta^{S}}(\mathcal{A}))=6>5.95=r_{3}^{\Delta^{S}}(\mathcal{A})r_{2}^{\overline{\Delta^{{S}}}}(\mathcal{A}).\nonumber
\end{eqnarray}
This implies that $\mathcal{A}$ satisfies the conditions (i)-(iv) in Theorem 4.1, thus $\mathcal{A}$ is positive definite.
\section{A new upper bound for the spectral radius of nonnegative tensors}
On the basis of the results in Section 3, we establish a new upper bound for the spectral radius of nonnegative tensors in this section, and compare this bound with some known bounds derived in \cite{4,13,31}.
\newtheorem{thm5}[thm1]{Theorem}
\begin{thm5}\label{thm5}
Let $\mathcal{A}\in{\mathbb{R}}^{[m,n]}$ be a nonnegative tensor with $n\geq{2}$. And let $S$ be a nonempty proper subset of $N$. Then
\begin{eqnarray}
\rho(\mathcal{A})\leq\eta_{\max}(\mathcal{A})=\max\{\eta_{1}(\mathcal{A}),\eta_{2}(\mathcal{A}),\eta_{3}(\mathcal{A}),\eta_{4}(\mathcal{A})\},\nonumber
\end{eqnarray}
where
\begin{eqnarray}
\eta_{1}(\mathcal{A})=\max\limits_{i\in{S}}\{a_{i\cdots i}+r_{i}^{\overline{\Delta^{\bar{S}}}}(\mathcal{A})\},\quad \eta_{2}(\mathcal{A})=\max\limits_{i\in{\bar{S}}}\{a_{i\cdots i}+r_{i}^{\overline{\Delta^{S}}}(\mathcal{A})\},\nonumber
\end{eqnarray}
and
\begin{eqnarray}
&&\eta_{3}(\mathcal{A})=\max\limits_{i\in{S},j\in{\bar{S}}}\min\left\{\frac{1}{2}\left(a_{i\cdots i}+a_{j\cdots j}+r_{i}^{\overline{\Delta^{\bar{S}}}}(\mathcal{A})+r_{j}^{\Delta^{\bar{S}}}(\mathcal{A})+\Phi_{i,j}^{\frac{1}{2}}(\mathcal{A})\right),R_{i}(\mathcal{A})\right\},\nonumber\\
&&\eta_{4}(\mathcal{A})=\max\limits_{i\in{\bar{S}},j\in{S}}\min\left\{\frac{1}{2}\left(a_{i\cdots i}+a_{j\cdots j}+r_{i}^{\overline{\Delta^{S}}}(\mathcal{A})+r_{j}^{\Delta^{S}}(\mathcal{A})+\Pi_{i,j}^{\frac{1}{2}}(\mathcal{A})\right),R_{i}(\mathcal{A})\right\},\nonumber
\end{eqnarray}
with
\begin{eqnarray}
&&\Phi_{i,j}(\mathcal{A})=(a_{i\cdots i}-a_{j\cdots j}+r_{i}^{\overline{\Delta^{\bar{S}}}}(\mathcal{A})-r_{j}^{\Delta^{\bar{S}}}(\mathcal{A}))^2+4r_{i}^{\Delta^{\bar{S}}}(\mathcal{A})r_{j}^{\overline{\Delta^{\bar{S}}}}(\mathcal{A}),\nonumber\\
&&\Pi_{i,j}(\mathcal{A})=(a_{i\cdots i}-a_{j\cdots j}+r_{i}^{\overline{\Delta^{S}}}(\mathcal{A})-r_{j}^{\Delta^{S}}(\mathcal{A}))^2+4r_{i}^{\Delta^{S}}(\mathcal{A})r_{j}^{\overline{\Delta^{S}}}(\mathcal{A}).\nonumber
\end{eqnarray}
\end{thm5}
\noindent {\bf Proof.}
Since $\mathcal{A}$ is a nonnegative tensor, from Lemma 2.1, we see that $\rho(\mathcal{A})$ is an eigenvalue of $\mathcal{A}$, it follows from Theorem 3.1 that
\begin{eqnarray}
\rho(\mathcal{A})\in\Upsilon^{S}(\mathcal{A}):=\Big(\Upsilon_{i,j}^{S}(\mathcal{A})\Big)\bigcup
\Big(\Upsilon_{i,j}^{\bar{S}}(\mathcal{A})\Big),
\end{eqnarray}
where
\begin{eqnarray}
&&\Upsilon_{i,j}^{S}(\mathcal{A})=\left(\bigcup_{i\in{S}}\hat{\Upsilon}_{i}^{1}(\mathcal{A})\right)\bigcup\left(\bigcup_{i\in{S},j\in{\bar{S}}}\Big(\tilde{\Upsilon}_{i,j}^{1}(\mathcal{A})\bigcap\Gamma_{i}(\mathcal{A})\Big)\right),\nonumber\\
&&\Upsilon_{i,j}^{\bar{S}}(\mathcal{A})=\left(\bigcup_{i\in{\bar{S}}}\hat{\Upsilon}_{i}^{2}(\mathcal{A})\right)\bigcup\left(\bigcup_{i\in{\bar{S}},j\in{S}}\Big(\tilde{\Upsilon}_{i,j}^{2}(\mathcal{A})\bigcap\Gamma_{i}(\mathcal{A})\Big)\right).\nonumber
\end{eqnarray}
If $\rho(\mathcal{A})\in\bigcup\limits_{i\in{S}}\hat{\Upsilon}_{i}^{1}(\mathcal{A})$, then there exists $i_{0}\in{S}$ such that $|\rho(\mathcal{A})-a_{i_{0}\cdots i_{0}}|\leq r_{i_{0}}^{\overline{\Delta^{\bar{S}}}}(\mathcal{A})$. Moreover, by Lemma 2.2, we see that $\rho(\mathcal{A})\geq\max\limits_{i\in{N}}\{a_{i\cdots i}\}$, then
\begin{eqnarray}
\rho(\mathcal{A})-a_{i_{0}\cdots i_{0}}\leq r_{i_{0}}^{\overline{\Delta^{\bar{S}}}}(\mathcal{A}),\nonumber
\end{eqnarray}
i.e.,
\begin{eqnarray}
\rho(\mathcal{A})\leq a_{i_{0}\cdots i_{0}}+r_{i_{0}}^{\overline{\Delta^{\bar{S}}}}(\mathcal{A})\leq\max\limits_{i\in{S}}\{a_{i\cdots i}+r_{i}^{\overline{\Delta^{\bar{S}}}}(\mathcal{A})\}.
\end{eqnarray}
If $\rho(\mathcal{A})\in\bigcup\limits_{i\in{\bar{S}}}\hat{\Upsilon}_{i}^{2}(\mathcal{A})$, then there is one $i_{1}\in{\bar{S}}$ such that $|\rho(\mathcal{A})-a_{i_{1}\cdots i_{1}}|\leq r_{i_{1}}^{\overline{\Delta^{S}}}(\mathcal{A})$, which implies that
\begin{eqnarray}
\rho(\mathcal{A})-a_{i_{1}\cdots i_{1}}\leq r_{i_{1}}^{\overline{\Delta^{S}}}(\mathcal{A}),\nonumber
\end{eqnarray}
i.e.,
\begin{eqnarray}
\rho(\mathcal{A})\leq a_{i_{1}\cdots i_{1}}+r_{i_{1}}^{\overline{\Delta^{S}}}(\mathcal{A})\leq\max\limits_{i\in{\bar{S}}}\{a_{i\cdots i}+r_{i}^{\overline{\Delta^{S}}}(\mathcal{A})\}.
\end{eqnarray}
For the case that $\rho(\mathcal{A})\in\bigcup\limits_{i\in{S},j\in{\bar{S}}}\Big(\tilde{\Upsilon}_{i,j}^{1}(\mathcal{A})\bigcap\Gamma_{i}(\mathcal{A})\Big)$, then there exist $p\in{S}$ and $q\in{\bar{S}}$ such that
\begin{eqnarray}
|\rho(\mathcal{A})-a_{p\cdots p}|\leq r_{p}(\mathcal{A})
\end{eqnarray}
and
\begin{eqnarray}
(|\rho(\mathcal{A})-a_{p\cdots p}|-r_{p}^{\overline{\Delta^{\bar{S}}}}(\mathcal{A}))(|\rho(\mathcal{A})-a_{q\cdots q}|-r_{q}^{\Delta^{\bar{S}}}(\mathcal{A}))
\leq r_{p}^{\Delta^{\bar{S}}}(\mathcal{A})r_{q}^{\overline{\Delta^{\bar{S}}}}(\mathcal{A}).
\end{eqnarray}
Combining Lemma 2.2 and (30) results in
\begin{eqnarray}
\rho(\mathcal{A})\leq a_{p\cdots p}+r_{p}(\mathcal{A})=R_{p}(\mathcal{A}).
\end{eqnarray}
Besides, by Lemma 2.2, we solve the quadratic Inequality (31) yields
\begin{eqnarray}
\rho(\mathcal{A})\leq\frac{1}{2}\{a_{p\cdots p}+a_{q\cdots q}+r_{p}^{\overline{\Delta^{\bar{S}}}}(\mathcal{A})+r_{q}^{\Delta^{\bar{S}}}(\mathcal{A})+\Phi_{p,q}^{\frac{1}{2}}(\mathcal{A})\},
\end{eqnarray}
where $\Phi_{p,q}(\mathcal{A})=(a_{p\cdots p}-a_{q\cdots q}+r_{p}^{\overline{\Delta^{\bar{S}}}}(\mathcal{A})-r_{q}^{\Delta^{\bar{S}}}(\mathcal{A}))^2+4r_{p}^{\Delta^{\bar{S}}}(\mathcal{A})r_{q}^{\overline{\Delta^{\bar{S}}}}(\mathcal{A})$.

Combining (32) and (33) results in
\setlength\arraycolsep{2pt}
\begin{eqnarray}
\rho(\mathcal{A})&\leq&\min\left\{\frac{1}{2}\left(a_{p\cdots p}+a_{q\cdots q}+r_{p}^{\overline{\Delta^{\bar{S}}}}(\mathcal{A})+r_{q}^{\Delta^{\bar{S}}}(\mathcal{A})+\Phi_{p,q}^{\frac{1}{2}}(\mathcal{A})\right),R_{p}(\mathcal{A})\right\}\nonumber\\
&\leq&\max\limits_{i\in{S},j\in{\bar{S}}}\min\left\{\frac{1}{2}\left(a_{i\cdots i}+a_{j\cdots j}+r_{i}^{\overline{\Delta^{\bar{S}}}}(\mathcal{A})+r_{j}^{\Delta^{\bar{S}}}(\mathcal{A})+\Phi_{i,j}^{\frac{1}{2}}(\mathcal{A})\right),R_{i}(\mathcal{A})\right\}.\nonumber\\
\end{eqnarray}
Furthermore, if $\rho(\mathcal{A})\in\bigcup\limits_{i\in{\bar{S}},j\in{S}}\Big(\tilde{\Upsilon}_{i,j}^{2}(\mathcal{A})\bigcap\Gamma_{i}(\mathcal{A})\Big)$,
then there exist $k\in{\bar{S}}$ and $l\in{S}$ such that
\begin{eqnarray}
|\rho(\mathcal{A})-a_{k\cdots k}|\leq r_{k}(\mathcal{A})
\end{eqnarray}
and
\begin{eqnarray}
(|\rho(\mathcal{A})-a_{k\cdots k}|-r_{k}^{\overline{\Delta^{{S}}}}(\mathcal{A}))(|\rho(\mathcal{A})-a_{l\cdots l}|-r_{l}^{\Delta^{S}}(\mathcal{A}))
\leq r_{k}^{\Delta^{S}}(\mathcal{A})r_{l}^{\overline{\Delta^{{S}}}}(\mathcal{A}).
\end{eqnarray}
Combining Lemma 2.2 and (35) gives
\begin{eqnarray}
\rho(\mathcal{A})\leq a_{k\cdots k}+r_{k}(\mathcal{A})=R_{k}(\mathcal{A}).
\end{eqnarray}
By Lemma 2.2,  Inequality (36) is equivalent to
\begin{eqnarray}
\rho(\mathcal{A})\leq\frac{1}{2}\{a_{k\cdots k}+a_{l\cdots l}+r_{k}^{\overline{\Delta^{S}}}(\mathcal{A})+r_{l}^{\Delta^{S}}(\mathcal{A})+\Pi_{k,l}^{\frac{1}{2}}(\mathcal{A})\},
\end{eqnarray}
where $\Pi_{k,l}(\mathcal{A})=(a_{k\cdots k}-a_{l\cdots l}+r_{k}^{\overline{\Delta^{S}}}(\mathcal{A})-r_{l}^{\Delta^{S}}(\mathcal{A}))^2+4r_{k}^{\Delta^{S}}(\mathcal{A})r_{l}^{\overline{\Delta^{S}}}(\mathcal{A})$.
By (37) and (38), we obtain
\setlength\arraycolsep{2pt}
\begin{eqnarray}
\rho(\mathcal{A})&\leq&\min\left\{\frac{1}{2}\left(a_{k\cdots k}+a_{l\cdots l}+r_{k}^{\overline{\Delta^{S}}}(\mathcal{A})+r_{l}^{\Delta^{S}}(\mathcal{A})+\Pi_{k,l}^{\frac{1}{2}}(\mathcal{A})\right),R_{k}(\mathcal{A})\right\}\nonumber\\
&\leq&\max\limits_{i\in{\bar{S}},j\in{S}}\min\left\{\frac{1}{2}\left(a_{i\cdots i}+a_{j\cdots j}+r_{i}^{\overline{\Delta^{S}}}(\mathcal{A})+r_{j}^{\Delta^{S}}(\mathcal{A})+\Pi_{i,j}^{\frac{1}{2}}(\mathcal{A})\right),R_{i}(\mathcal{A})\right\}.\nonumber\\
\end{eqnarray}
The conclusion follows from Inequalities (28), (29), (34) and (39). \hfill$\blacksquare$
\newtheorem{rem6}[rem1]{Remark}
\begin{rem6}\label{rem4}
\emph{As the upper bounds for $\rho(\mathcal{A})$ in Theorems 3.3 and 3.4 in \cite{31} deduced from the eigenvalue localization sets $\mathcal{K}(\mathcal{A})$ and $\mathcal{K}^{S}(\mathcal{A})$, respectively, and that in Theorem 5.1 derived from the eigenvalue localization set $\Upsilon^{S}(\mathcal{A})$. It follows from $\Upsilon^{S}(\mathcal{A})\subseteq\Omega^{S}(\mathcal{A})\subseteq\mathcal{K}^{S}(\mathcal{A})\subseteq\mathcal{K}(\mathcal{A})\subseteq\Gamma(\mathcal{A})$ and the fact that $\omega^{S}_{\max}(\mathcal{A})\leq\omega_{\max}(\mathcal{A})\leq {R_{\max}(\mathcal{A})}$ (see Theorems 3.5 in \cite{31}) that $\eta_{\max}(\mathcal{A})\leq\omega^{S}_{\max}(\mathcal{A})\leq\omega_{\max}(\mathcal{A})\leq {R_{\max}(\mathcal{A})}$, we use $\omega_{\max}(\mathcal{A})$ and $\omega^{S}_{\max}(\mathcal{A})$ to denote the upper bounds in Theorems 3.3 and 3.4 in \cite{31}, respectively in this paper.  }
\end{rem6}

For some upper bounds we have showed that our bound is sharper than existing bounds. Now we take an example to show the efficiency of the new upper bounds.
\begin{Exa}\label{E2}
\setlength\arraycolsep{4pt}
\emph{Consider the following nonnegative tensor}
\begin{eqnarray}
\mathcal{A}=[A(1,:,:),A(2,:,:),A(3,:,:)]\in{\mathbb{R}}^{[3,3]},\nonumber
\end{eqnarray}
\emph{where}
\setlength\arraycolsep{4pt}
\begin{eqnarray}
A(1,:,:)=\left(
      \begin{array}{ccc}
        3 & 1 & 0 \\
        0 & 1 & 2 \\
        0 & 0 & 2 \\
      \end{array}
    \right),
A(2,:,:)=\left(
      \begin{array}{ccc}
        2 & 0 & 3 \\
        0 & 1 & 0 \\
        0 & 0 & 1 \\
      \end{array}
    \right),
A(3,:,:)=\left(
           \begin{array}{ccc}
             15 & 1 & 8 \\
             4 & 1 & 0 \\
             0 & 0 & 1 \\
           \end{array}
         \right).\nonumber
\end{eqnarray}
\end{Exa}
We compare the results derived in Theorem 5.1 with those in Lemma 5.2 of \cite{4}, Theorems 3.3 and 3.4 of \cite{31}. Let $S=\{1,2\}$, then $\bar{S}=\{3\}$. By Lemma 5.2 of \cite{4}, we have
$$\rho(\mathcal{A})\leq 30.$$
By Theorems 3.3 and 3.4 of \cite{31}, we have
$$\rho(\mathcal{A})\leq 29.2127.$$
By Theorem 13 of \cite{13}, we get
$$\rho(\mathcal{A})\leq 20.2250.$$
By Theorem 5.1, we obtain
$$\rho(\mathcal{A})\leq 15.6437.$$
This shows that the upper bound in Theorem 5.1 is sharper than those in Lemma 5.2 of \cite{4} and Theorems 3.3-3.4 of \cite{31}, and better than the one in Theorem 13 of \cite{13} in some cases.
\section{A new lower bound for the minimum \emph{H}-eigenvalue of weakly irreducible strong \emph{M}-tensors}
In this section, by applying the results of Theorem 3.1, we exhibit a new lower bound for the minimum \emph{H}-eigenvalue of weakly irreducible strong \emph{M}-tensors, which improves some existing ones derived in \cite{6,7}.
\newtheorem{thm4}[thm1]{Theorem}
\begin{thm4}\label{thm4}
Let $\mathcal{A}\in{\mathbb{R}}^{[m,n]}$ be a weakly irreducible strong M-tensor with $n\geq{2}$. And let $S$ be a nonempty proper subset of $N$. Then
\begin{eqnarray}
\tau(\mathcal{A})\geq\pi_{\min}(\mathcal{A})=\min\{\pi_{1}(\mathcal{A}),\pi_{2}(\mathcal{A}),\pi_{3}(\mathcal{A}),\pi_{4}(\mathcal{A})\},\nonumber
\end{eqnarray}
where
\begin{eqnarray}
\pi_{1}(\mathcal{A})=\min\limits_{i\in{S}}\{a_{i\cdots i}-r_{i}^{\overline{\Delta^{\bar{S}}}}(\mathcal{A})\},\quad \pi_{2}(\mathcal{A})=\min\limits_{i\in{\bar{S}}}\{a_{i\cdots i}-r_{i}^{\overline{\Delta^{S}}}(\mathcal{A})\},\nonumber
\end{eqnarray}
and
\begin{eqnarray}
&&\pi_{3}(\mathcal{A})=\min\limits_{i\in{S},j\in{\bar{S}}}\max\left\{\frac{1}{2}\left(a_{i\cdots i}+a_{j\cdots j}-r_{i}^{\overline{\Delta^{\bar{S}}}}(\mathcal{A})-r_{j}^{\Delta^{\bar{S}}}(\mathcal{A})-\Theta_{i,j}^{\frac{1}{2}}(\mathcal{A})\right),R_{i}(\mathcal{A})\right\},\nonumber\\
&&\pi_{4}(\mathcal{A})=\min\limits_{i\in{\bar{S}},j\in{S}}\max\left\{\frac{1}{2}\left(a_{i\cdots i}+a_{j\cdots j}-r_{i}^{\overline{\Delta^{S}}}(\mathcal{A})-r_{j}^{\Delta^{S}}(\mathcal{A})-\Lambda_{i,j}^{\frac{1}{2}}(\mathcal{A})\right),R_{i}(\mathcal{A})\right\},\nonumber
\end{eqnarray}
with
\begin{eqnarray}
&&\Theta_{i,j}(\mathcal{A})=(a_{i\cdots i}-a_{j\cdots j}-r_{i}^{\overline{\Delta^{\bar{S}}}}(\mathcal{A})+r_{j}^{\Delta^{\bar{S}}}(\mathcal{A}))^2+4r_{i}^{\Delta^{\bar{S}}}(\mathcal{A})r_{j}^{\overline{\Delta^{\bar{S}}}}(\mathcal{A}),\nonumber\\
&&\Lambda_{i,j}(\mathcal{A})=(a_{i\cdots i}-a_{j\cdots j}-r_{i}^{\overline{\Delta^{S}}}(\mathcal{A})+r_{j}^{\Delta^{S}}(\mathcal{A}))^2+4r_{i}^{\Delta^{S}}(\mathcal{A})r_{j}^{\overline{\Delta^{S}}}(\mathcal{A}).\nonumber
\end{eqnarray}
\end{thm4}
\noindent {\bf Proof.}
Inasmuch as $\mathcal{A}$ is a weakly irreducible strong \emph{M}-tensor, from Lemma 2.3, $\tau(\mathcal{A})$ is an eigenvalue of $\mathcal{A}$, by Theorem 3.1, we have
\begin{eqnarray}
\tau(\mathcal{A})\in\Upsilon^{S}(\mathcal{A}):=\Big(\Upsilon_{i,j}^{S}(\mathcal{A})\Big)\bigcup
\Big(\Upsilon_{i,j}^{\bar{S}}(\mathcal{A})\Big),\nonumber
\end{eqnarray}
where
\begin{eqnarray}
&&\Upsilon_{i,j}^{S}(\mathcal{A})=\left(\bigcup_{i\in{S}}\hat{\Upsilon}_{i}^{1}(\mathcal{A})\right)\bigcup\left(\bigcup_{i\in{S},j\in{\bar{S}}}\Big(\tilde{\Upsilon}_{i,j}^{1}(\mathcal{A})\bigcap\Gamma_{i}(\mathcal{A})\Big)\right),\nonumber\\
&&\Upsilon_{i,j}^{\bar{S}}(\mathcal{A})=\left(\bigcup_{i\in{\bar{S}}}\hat{\Upsilon}_{i}^{2}(\mathcal{A})\right)\bigcup\left(\bigcup_{i\in{\bar{S}},j\in{S}}\Big(\tilde{\Upsilon}_{i,j}^{2}(\mathcal{A})\bigcap\Gamma_{i}(\mathcal{A})\Big)\right).\nonumber
\end{eqnarray}
If $\tau(\mathcal{A})\in\bigcup\limits_{i\in{S}}\hat{\Upsilon}_{i}^{1}(\mathcal{A})$, then there exists $i_{0}\in{S}$ such that $|\tau(\mathcal{A})-a_{i_{0}\cdots i_{0}}|\leq r_{i_{0}}^{\overline{\Delta^{\bar{S}}}}(\mathcal{A})$. Besides, using Lemma 2.4, we know that $\tau(\mathcal{A})\leq\min\limits_{i\in{N}}\{a_{i\cdots i}\}$, then
\begin{eqnarray}
a_{i_{0}\cdots i_{0}}-\tau(\mathcal{A})\leq r_{i_{0}}^{\overline{\Delta^{\bar{S}}}}(\mathcal{A}),\nonumber
\end{eqnarray}
i.e.,
\begin{eqnarray}
\tau(\mathcal{A})\geq a_{i_{0}\cdots i_{0}}-r_{i_{0}}^{\overline{\Delta^{\bar{S}}}}(\mathcal{A})\geq\min\limits_{i\in{S}}\{a_{i\cdots i}-r_{i}^{\overline{\Delta^{\bar{S}}}}(\mathcal{A})\}.
\end{eqnarray}
If $\tau(\mathcal{A})\in\bigcup\limits_{i\in{\bar{S}}}\hat{\Upsilon}_{i}^{2}(\mathcal{A})$, then there is one $i_{1}\in{\bar{S}}$ such that $|\tau(\mathcal{A})-a_{i_{1}\cdots i_{1}}|\leq r_{i_{1}}^{\overline{\Delta^{S}}}(\mathcal{A})$, together with Lemma 2.4 yields
\begin{eqnarray}
a_{i_{1}\cdots i_{1}}-\tau(\mathcal{A})\leq r_{i_{1}}^{\overline{\Delta^{S}}}(\mathcal{A}),\nonumber
\end{eqnarray}
i.e.,
\begin{eqnarray}
\tau(\mathcal{A})\geq a_{i_{1}\cdots i_{1}}-r_{i_{1}}^{\overline{\Delta^{S}}}(\mathcal{A})\geq\min\limits_{i\in{\bar{S}}}\{a_{i\cdots i}-r_{i}^{\overline{\Delta^{S}}}(\mathcal{A})\}.
\end{eqnarray}
For the case that $\tau(\mathcal{A})\in\bigcup\limits_{i\in{S},j\in{\bar{S}}}\Big(\tilde{\Upsilon}_{i,j}^{1}(\mathcal{A})\bigcap\Gamma_{i}(\mathcal{A})\Big)$, then there exist $p\in{S}$ and $q\in{\bar{S}}$ such that
\begin{eqnarray}
|\tau(\mathcal{A})-a_{p\cdots p}|\leq r_{p}(\mathcal{A})
\end{eqnarray}
and
\begin{eqnarray}
(|\tau(\mathcal{A})-a_{p\cdots p}|-r_{p}^{\overline{\Delta^{\bar{S}}}}(\mathcal{A}))(|\tau(\mathcal{A})-a_{q\cdots q}|-r_{q}^{\Delta^{\bar{S}}}(\mathcal{A}))
\leq r_{p}^{\Delta^{\bar{S}}}(\mathcal{A})r_{q}^{\overline{\Delta^{\bar{S}}}}(\mathcal{A}).
\end{eqnarray}
Combining Lemma 2.4 and (42) gives
\begin{eqnarray}
\tau(\mathcal{A})\geq a_{p\cdots p}-r_{p}(\mathcal{A})=R_{p}(\mathcal{A}).
\end{eqnarray}
Having in mind that $\tau(\mathcal{A})\leq\min\limits_{i\in{N}}\{a_{i\cdots i}\}$, it follows from (43) that
\begin{eqnarray}
\tau(\mathcal{A})\geq\frac{1}{2}\{a_{p\cdots p}+a_{q\cdots q}-r_{p}^{\overline{\Delta^{\bar{S}}}}(\mathcal{A})-r_{q}^{\Delta^{\bar{S}}}(\mathcal{A})-\Theta_{p,q}^{\frac{1}{2}}(\mathcal{A})\},
\end{eqnarray}
where $\Theta_{p,q}(\mathcal{A})=(a_{p\cdots p}-a_{q\cdots q}-r_{p}^{\overline{\Delta^{\bar{S}}}}(\mathcal{A})+r_{q}^{\Delta^{\bar{S}}}(\mathcal{A}))^2+4r_{p}^{\Delta^{\bar{S}}}(\mathcal{A})r_{q}^{\overline{\Delta^{\bar{S}}}}(\mathcal{A})$.

Combining (44) and (45) results in
\setlength\arraycolsep{2pt}
\begin{eqnarray}
\tau(\mathcal{A})&\geq&\max\left\{\frac{1}{2}\left(a_{p\cdots p}+a_{q\cdots q}-r_{p}^{\overline{\Delta^{\bar{S}}}}(\mathcal{A})-r_{q}^{\Delta^{\bar{S}}}(\mathcal{A})-\Theta_{p,q}^{\frac{1}{2}}(\mathcal{A})\right),R_{p}(\mathcal{A})\right\}\nonumber\\
&\geq&\min\limits_{i\in{S},j\in{\bar{S}}}\max\left\{\frac{1}{2}\left(a_{i\cdots i}+a_{j\cdots j}-r_{i}^{\overline{\Delta^{\bar{S}}}}(\mathcal{A})-r_{j}^{\Delta^{\bar{S}}}(\mathcal{A})-\Theta_{i,j}^{\frac{1}{2}}(\mathcal{A})\right),R_{i}(\mathcal{A})\right\}.\nonumber\\
\end{eqnarray}
Furthermore, if $\tau(\mathcal{A})\in\bigcup\limits_{i\in{\bar{S}},j\in{S}}\Big(\tilde{\Upsilon}_{i,j}^{2}(\mathcal{A})\bigcap\Gamma_{i}(\mathcal{A})\Big)$,
then there exist $k\in{\bar{S}}$ and $l\in{S}$ such that
\begin{eqnarray}
|\tau(\mathcal{A})-a_{k\cdots k}|\leq r_{k}(\mathcal{A})
\end{eqnarray}
and
\begin{eqnarray}
(|\tau(\mathcal{A})-a_{k\cdots k}|-r_{k}^{\overline{\Delta^{{S}}}}(\mathcal{A}))(|\tau(\mathcal{A})-a_{l\cdots l}|-r_{l}^{\Delta^{S}}(\mathcal{A}))
\leq r_{k}^{\Delta^{S}}(\mathcal{A})r_{l}^{\overline{\Delta^{{S}}}}(\mathcal{A}).
\end{eqnarray}
It follows from Lemma 2.4 and (47) that
\begin{eqnarray}
\tau(\mathcal{A})\geq a_{k\cdots k}-r_{k}(\mathcal{A})=R_{k}(\mathcal{A}).
\end{eqnarray}
On the other hand, by Lemma 2.4, solving $\tau(\mathcal{A})$ in Inequality (48) yields
\begin{eqnarray}
\tau(\mathcal{A})\geq\frac{1}{2}\{a_{k\cdots k}+a_{l\cdots l}-r_{k}^{\overline{\Delta^{S}}}(\mathcal{A})-r_{l}^{\Delta^{S}}(\mathcal{A})-\Lambda_{k,l}^{\frac{1}{2}}(\mathcal{A})\},
\end{eqnarray}
where $\Lambda_{k,l}(\mathcal{A})=(a_{k\cdots k}-a_{l\cdots l}-r_{k}^{\overline{\Delta^{S}}}(\mathcal{A})+r_{l}^{\Delta^{S}}(\mathcal{A}))^2+4r_{k}^{\Delta^{S}}(\mathcal{A})r_{l}^{\overline{\Delta^{S}}}(\mathcal{A})$,
which together with (49) gives
\setlength\arraycolsep{2pt}
\begin{eqnarray}
\tau(\mathcal{A})&\geq&\max\left\{\frac{1}{2}\left(a_{k\cdots k}+a_{l\cdots l}-r_{k}^{\overline{\Delta^{S}}}(\mathcal{A})-r_{l}^{\Delta^{S}}(\mathcal{A})-\Lambda_{k,l}^{\frac{1}{2}}(\mathcal{A})\right),R_{k}(\mathcal{A})\right\}\nonumber\\
&\geq&\min\limits_{i\in{\bar{S}},j\in{S}}\max\left\{\frac{1}{2}\left(a_{i\cdots i}+a_{j\cdots j}-r_{i}^{\overline{\Delta^{S}}}(\mathcal{A})-r_{j}^{\Delta^{S}}(\mathcal{A})-\Lambda_{i,j}^{\frac{1}{2}}(\mathcal{A})\right),R_{i}(\mathcal{A})\right\}.\nonumber\\
\end{eqnarray}
The results of this theorem follow from the Inequalities (40), (41), (46) and (51). This proves the theorem. \hfill$\blacksquare$

\newtheorem{rem5}[rem1]{Remark}
\begin{rem5}\label{rem4}
\emph{Inasmuch as the lower bounds for $\tau(\mathcal{A})$ in Theorem 2.2 in \cite{6} and Theorem 6.1 derived from the eigenvalue localization sets $\mathcal{K}(\mathcal{A})$ and $\Upsilon^{S}(\mathcal{A})$, respectively. Using the similar technique as Theorem 3.5 of \cite{31}, we can prove that $\varrho_{\min}(\mathcal{A})\geq {R_{\min}(\mathcal{A})}$, here we use $\varrho_{\min}(\mathcal{A})$ to denote the lower bound in Theorem 2.2 in \cite{6}. Combining $\Upsilon^{S}(\mathcal{A})\subseteq\Omega^{S}(\mathcal{A})\subseteq\mathcal{K}^{S}(\mathcal{A})\subseteq\mathcal{K}(\mathcal{A})\subseteq\Gamma(\mathcal{A})$ with $\varrho_{\min}(\mathcal{A})\geq {R_{\min}(\mathcal{A})}$ results in $\pi_{\min}(\mathcal{A})\geq\varrho_{\min}(\mathcal{A})\geq {R_{\min}(\mathcal{A})}$, i.e., the lower bound in Theorem 6.1 is an improvement on those in Theorems 2.1-2.2 of \cite{6}. }
\end{rem5}
Let us show that by a simple example as follows.
\begin{Exa}\label{E2}
\setlength\arraycolsep{4pt}
\emph{Consider the following irreducible nonsingular \emph{M}-tensor}
\begin{eqnarray}
\mathcal{A}=[A(1,:,:),A(2,:,:),A(3,:,:)]\in{\mathbb{R}}^{[3,3]},\nonumber
\end{eqnarray}
\emph{where}
\setlength\arraycolsep{4pt}
\begin{eqnarray}
&&A(1,:,:)=\left(
      \begin{array}{ccc}
        12 & -2.2 & -0.3 \\
        0 & 0 & -2 \\
        0 & -1 & -1.5 \\
      \end{array}
    \right),\nonumber\\
&&A(2,:,:)=\left(
      \begin{array}{ccc}
        -0.5 & -4.8 & -8 \\
        0 & 30 & 0 \\
        -1 & 0 & -0.5 \\
      \end{array}
    \right),\nonumber\\
&&A(3,:,:)=\left(
           \begin{array}{ccc}
             0 & -3 & -1 \\
             0 & -1 & -3.5 \\
             -1 & -3 & 15 \\
           \end{array}
         \right).\nonumber
\end{eqnarray}
\end{Exa}
We compare the results exhibited in Theorem 6.1 with those in Theorems 2.1-2.2 of \cite{6} and Theorem 4.5 of \cite{7}. Let $S=\{1,2\}$, then $\bar{S}=\{3\}$. By Theorems 2.1-2.2 of \cite{6}, we have
$$\tau(\mathcal{A})\geq 2.5.$$
By Theorem 4.5 of \cite{7}, we get
$$\tau(\mathcal{A})\geq 2.74.$$
By Theorem 6.1, we obtain
$$\tau(\mathcal{A})\geq 6.5,$$
which shows that the lower bound in Theorem 6.1 is much better than those in Theorems 2.1-2.2 of \cite{6} and Theorem 4.5 of \cite{7}.
\section{Concluding remarks}
In this paper, a new \emph{S}-type eigenvalue localization set for tensors is established,
which is proved to be sharper than the ones in \cite{16,31}. As applications of this new set, checkable sufficient conditions for the positive definiteness and the positive semi-definiteness of tensors are proposed, these conditions have wider scope of applications compare with those of \cite{16,21,31}. Moreover, based on the results of Theorem 3.1, we give
new bounds for the spectral radius of nonnegative tensors and the minimum
\emph{H}-eigenvalue of weakly irreducible strong \emph{M}-tensors, these bounds improve some existing ones
obtained by Yang and Yang \cite{4}, Li et al. \cite{31} and He and Huang \cite{6}. Numerical experiments are also implemented to illustrate the advantages of these results.

However, the new \emph{S}-type eigenvalue localization set and the derived bounds
depend on the set \emph{S}. How to choose \emph{S} to make $\Upsilon^{S}(\mathcal{A})$ and the bounds exhibited
in this paper as tight as possible is very important and interesting, while if the
dimension of the tensor $\mathcal{A}$ is large, this work is very difficult. Therefore, future
work will include numerical or theoretical studies for finding the best choice for
\emph{S}.



\end{document}